\documentclass[a4paper,twoside,
11pt]{article}
\usepackage{amssymb,amsfonts,amsmath}
\begin{document}
\newtheorem{defin}{~~~~Definition}
\newtheorem{prop}{~~~~Proposition}
\newtheorem{remark}{~~~~Remark}
\newtheorem{cor}{~~~~Corollary}
\newtheorem{theor}{~~~~Theorem}
\newtheorem{lemma}{~~~~Lemma}
\newtheorem{ass}{~~~~Assumption}
\newtheorem{con}{~~~~Conjecture}
\newtheorem{concl}{~~~~Conclusion}
\renewcommand{\theequation}{\thesection.\arabic{equation}}

\title{
On geodesic equivalence of Riemannian metrics and 
sub-Riemannian metrics on distributions of corank 1} 
\date{}
\author{Igor Zelenko\thanks{S.I.S.S.A., Via Beirut 2-4,
34014, Trieste, Italy; email: zelenko@sissa.it}} \maketitle
\setcounter{equation}{0}

\begin{abstract}
The present paper is devoted to the problem of (local)
geodesic equivalence of Riemannian metrics and
sub-Riemannian metrics on generic corank 1 distributions.
Using Pontryagin Maximum Principle, we treat Riemannian and 
sub-Riemannian cases in an unified way and obtain some 
algebraic necessary conditions for the geodesic equivalence 
of (sub-)Riemannian metrics. In this way first we obtain a 
new elementary proof of classical Levi-Civita's Theorem 
about the classification of all Riemannian geodesically 
equivalent metrics in a neighborhood of so-called regular 
(stable) point w.r.t. these metrics. Secondly we prove that 
sub-Riemannian metrics on contact 
distributions are geodesically equivalent iff they are 
constantly proportional. Then we describe all geodesically 
equivalent sub-Riemannian metrics on quasi-contact 
distributions. 
Finally we make the classification of all pairs of
geodesically equivalent Riemannian metrics on a surface,
which proportional in an isolated point. This is the
simplest case, which was not covered by Levi-Civita's
Theorem.
\end{abstract}

\section{Introduction}
\setcounter{equation}{0} \indent

Let us recall that two Riemannian metrics on a
manifold $M$ are called {\it geodesically} (or {\it
projective}) {\it equivalent}
 at a point $q_0\in M$, if
in some neighborhood of $q_0$ all their geodesics,
considered as unparametrized curves, coincide.
 The notion of geodesic equivalence can be generalized
directly to sub-Riemannian metrics by replacing Riemannian
geodesics by normal sub-Riemannian geodesics:

Let $D$ be a bracket-generating (completely nonholonomic) 
distribution on $M$. A Lipschitzian curve $\xi(t) $ is 
called admissible for the distribution $D$, if it is 
tangent to $D$ almost everywhere, i.e., $\dot \xi(t)\in 
D\bigl(\xi(t)\bigr)$ a.e.. 
 A sub-Riemannian metric $G$ on
$D$ is given by choosing an inner product
$G_q(\cdot,\cdot)$ on each subspaces $D(q)$ for any $q\in
M$ smoothly w.r.t. $q$. Let 
$||\cdot||_q=\sqrt{G_q(\cdot,\cdot)}$ be the corresponding 
Euclidean norm on $D(q)$. For any admissible curve 
$\xi:[0,T]\mapsto M$ its length w.r.t. the sub-Riemannian 
metric $G$ is equal to $\int_0^T 
||\dot\xi(t)||_{\xi(t)} \,d t$. Given two points $q_1$ and
$q_2$ one can look for the curve of minimal length among
all admissible curves connecting $q_1$ with $q_2$. This
problem can be obviously reformulated as a time-minimal
control problem (for this one takes into the consideration 
only admissible curves parametrized by the length). The 
{\it sub-Riemannian extremal trajectory} w.r.t. the metric 
$G$ is the projection to $M$ of a Pontryagin extremal of 
this problem (which lives in the cotangent bundle $T^*M$). 

In general, Pontryagin extremals can be normal or abnormal:
the extremal is called abnormal, if the Lagrange multiplier
of the functional is equal to zero, and normal otherwise.
The projection of normal (abnormal) Pontryagin extremal is
called a {\it normal (abnormal) sub-Riemannian extremal
trajectory}. Any abnormal sub-Riemannian extremal
trajectory, considered as unparametrized curves, is
characterized by distribution $D$ only, but not by the
metric on it. Normal sub-Riemannian extremals surely depend 
on the metric. They can be described in the following 
simple way: Let $h:T^*M\mapsto R$ satisfies 
\begin{equation}
\label{hamsR}
 h(p,q)=\frac{1}{2}\bigl(\max\{p(v):||v||_q=1, v\in
 D(q)\}\bigr)^2 \quad q\in M,\, p\in
 T_q^*M.
 \end{equation}
Then the normal sub-Riemannian extremal trajectories are
exactly the projections on $M$ of the trajectories of the
Hamiltonian system $\dot\lambda=\vec h(\lambda)$, lying on
the $\frac{1}{2}$-level set of $h$, i.e., on the set
$\{\lambda\in T^*M
:h(\lambda)=\frac{1}{2}\}$.

 \begin{remark}
 The norm $||\cdot||_q$ on $D_q$ induces the norm on the
 dual space, which will be denoted also by $||\cdot||_q$.
 Therefore taking the restriction $p|_{D(q)}$ of some covector $p\in T_q^*M$ one can rewrite (\ref{hamsR})
 in the  following form
 \begin{equation}
\label{hamsR1}
 h(p,q)=\frac{1}{2}||p|_{D(q)}||_q^2\quad q\in M,\, p\in
 T_q^*M.
 \end{equation}
 \end{remark}

Note that a Riemannian metric is actually the 
sub-Riemannian metric with $D=TM$ and classical Riemannian 
geodesics are exactly normal extremal trajectories in this 
situation (here abnormal extremals do not exist). Note also 
that, as in Riemannian case, sufficiently small pieces of 
normal sub-Riemannian extremal trajectories are length 
minimizers (see, for example, \cite{sus}, Appendix C 
there). Therefore we will call them in the sequel {\it 
normal sub-Riemannian geodesics}. The following definition 
is a natural extension of the notion of the geodesic 
equivalence from Riemannian to the general sub-Riemannian 
case: 

\begin{defin}
\label{subdef}
 Two sub-Riemannian metrics given on a distribution $D$
of a manifold $M$ are called {\it geodesically} (or {\it
projective}) {\it equivalent}
 at a point $q_0\in M$, if
in some neighborhood of $q_0$ all their normal geodesics,
considered as unparametrized curves, coincide.
\end{defin}

It is clear that if sub-Riemannian metrics $G_1$ and $G_2$
are constantly proportional, i.e., there exists a positive
constant $C$ such that ${G_2}_q=C {G_2}_q$ for any $q$,
then they are geodesically equivalent. The first appearing
question is whether there exist constantly non-proportional
geodesically equivalent sub-Riemannian metrics? The
simplest example of constantly non-proportional Riemannian
metrics on a surface can be described as follows: Let $P$
and $S$ be a plane and a hemisphere in $\mathbb R^3$ such
that equator of the hemisphere is parallel to the plane.
Let $G_1$ and $\bar G_2$ be the metrics on $P$ and $S$
respectively, induced from the Euclidean metric on $\mathbb
R^3$. Denote by $F:S\mapsto P$ the stereographic projection
from the center $O$ of the hemisphere (namely, if $q\in S$
then $F(q)$ is the only point on $P$ lying on the straight
line, which connects $O$ and $q$). Then the mapping $F$
sends geodesics of $\bar G_2$ (arcs of big circles on $S$)
to geodesics of $G_1$ (straight lines on $P$). Therefore
$G_1$ is geodesically equivalent to $G_2=(F^{-1})^*\bar
G_2$, the pull-back of $\bar G_2$ by $F^{-1}$, but this
metric are not constantly proportional. Moreover, as E.
Beltrami showed in \cite{belt}, a Riemannian metric on a
surface is geodesically equivalent to the flat one iff it
has a constant curvature.

Let us introduced some notions, which are important for the 
considered problem. For a given ordered pair of 
sub-Riemannian metrics $G_1$, $G_2$ and a point $q$ one can 
define the following linear operator $S_q:D(q)\mapsto 
D(q)$: $${G_2}_q(v_1,v_2)={G_1}_q(S_q v_1,v_2),\quad 
v_1,v_2\in D(q).$$ Obviously, $S_q$ is self-adjoint w.r.t. 
the Euclidean structure given by $G_1$. 
\begin{defin}
\label{transdef}
The operator $S_q$ will
be called the {\it transition operator} from the metric
$G_1$ to the metric $G_2$ at the point $q$.
\end{defin}

Let $N(q)$ be the number of distinct eigenvalues of the 
operator $S_q$. 

\begin{defin}
\label{regdef}
The point $q_0$ is called {\it regular} w.r.t. the pair of
sub-Riemannian metrics $G_1$ and $G_2$, if the function
$N(q)$ is constant in some neighborhood of $q_0$.
\end{defin}

Note that the regularity of the point $q_0$ is equivalent 
to the fact that the set of multiplicities of eigenvalues 
of the transition operator $S_q$ is the same for all points 
$q$ from some neighborhood of $q_0$ (in \cite{top} regular 
points were called stable). By standard arguments one can 
show that the function $N(q)$ is lower semicontinuous. This 
together with the fact that it is integer-valued implies 
the following 

\begin{prop}
\label{denseprop} The set of regular points w.r.t. the pair
of sub-Riemannian metrics is open and dense in $M$.
\end{prop}

For Riemannian metrics on an $n$-dimensional manifold all
possible pairs of geodesically equivalent metrics in a
neighborhood of a regular
point w.r.t. these metrics were described already by 
Levi-Civita in \cite{levi} (see Theorem \ref {dinit} below 
and also \cite{top}), who had extended the earlier result 
of Dini for surfaces (see \cite{dini},\cite{eisen}, or 
\cite{mat}) to an arbitrary $n$. From this result it 
follows that Riemannian metrics, for which there exists at 
least one non-proportional geodesically equivalent 
Riemannian metric, are of the very special form. 

The classification of geodesically equivalent Riemannian
metrics at non-regular points (i.e., points, where
eigenvalues of transition operator bifurcate) even on a
surface was not done, while the geodesic equivalence of
proper sub-Riemannian metrics (i.e, when $D\neq TM$ and $D$ 
is bracket-generating) was not studied before. In the 
present paper we treat both these problems.

In the sequel for shortness $(m,n)$-distribution means an
$m$-dimensional subbundle of the tangent bundle of an
$n$-dimensional manifold. Our study of the geodesic 
equivalence of proper sub-Riemannian metrics will be mainly 
concentrated on the following two cases: 
\begin{enumerate}
\item
$D$ is the contact distribution. Namely, $D$ is a corank 1
distribution on an odd dimensional manifold such that if
$\omega$ is a differential $1$-form, which annihilates $D$,
$D(q)=\{v\in T_q M:\omega_q(v)=0\}$, then the restriction
$d\omega|_D$ of the differential $d\omega$ on $D$ is a
nondegenerated $2$-form at any $q$. In this case there are
no abnormal Pontryagin extremals.

\item
$D$ is the quasi-contact distribution. Namely, $D$ is a
corank 1 distribution on an even dimensional manifold such
that if $\omega$ is a differential $1$-form, which
annihilates $D$, then the restriction $d\omega|_D$ of the
differential $d\omega$ on $D$ has $1$-dimensional kernel at
any $q$. The kernels of $d\omega|_D$ form line 
distribution. We will call it {\it the abnormal line 
distribution of the quasi-contact distribution} $D$. 
Abnormal extremal trajectories of the sub-Riemannian metric 
$G$ on $D$ are exactly the leaves of this distribution, 
parametrized by the length. 
\end{enumerate}

Clearly in both cases the germs of distribution $D$ are 
generic germs of corank 1 distributions. Note also that in 
both cases there exists only one, up to a diffeomorphism, 
distribution satisfying the prescribed properties (a 
particular case of Darboux's Theorem). Actually, our method 
works for sub-Riemannian metrics defined on much more 
general distributions, for example, on so-called step $1$ 
bracket-generating distributions: an $(m,n)$-distribution 
$D$ is called step $1$ bracket-generating if $\dim 
D^{l+1}=\dim D^{l}+1$ for any $1\leq l\leq n-m$ (here the 
$l$th power $D^l$ of the distribution $D$ is defined by 
induction $D^l=D^{l-1}+[D, D^{l-1}]$). 
 The study of the 
problem of the geodesic equivalence for sub-Riemannian 
metrics on general step $1$ bracket-generating 
distributions will be done in our future publications. 

The paper is organized as follows. In section 2 we show 
that the problem of geodesic equivalence of sub-Riemannian 
metrics can be reduced to the question of the existence of 
an orbital diffeomorphism between the corresponding flows 
of extremals. This reduction is obvious in the Riemannian 
case, but in the proper sub-Riemannian case it has some 
additional difficulties, especially in the presence of 
abnormal extremals. After this reduction we express the 
condition for the existence of the orbital diffeomorphism 
in terms of the special frame adapted to the pair of 
sub-Riemannian metrics. Further for step $1$ 
bracket-generating distributions we obtain a necessary 
condition for the geodesic equivalence in terms of 
divisibility of some polynomials (on the fibers of the 
cotangent bundle of the ambient manifold) associated with 
these metrics. We call it the first divisibility condition. 
It imposes rather 
strong restrictions on the pair of the metrics.

In section 3 we give the coordinate-free formulation of 
Levi-Civita's Theorem (Theorem \ref{dinit}) and prove it in 
a new, rather elementary way, using the conditions for the 
existence of the orbital diffeomorphism and the first 
divisibility condition. In section 4 for a sub-Riemannian 
metric on corank 1 distribution we obtain an additional 
necessary condition for the geodesic equivalence in terms 
of divisibility of some polynomials associated with these 
metrics. We call it the second divisibility condition. 
Using the conditions for the existence of the orbital 
diffeomorphism and the second divisibility condition we 
prove that sub-Riemannian metrics on contact 
distributions are geodesically equivalent iff they are 
constantly proportional (Theorem \ref{contprop}) and we 
give the classification of all geodesically equivalent 
sub-Riemannian metrics on quasi-contact distributions 
(Theorem \ref{quasit}). This classification is given in 
coordinate-free way and has apparent similarities with our 
interpretation of Levi-Civita's Theorem. This gives a hope 
for the existence of a general classification theorem about 
geodesic equivalence of sub-Riemannian metrics defined on 
very general class of distribution, which will contain as 
particular cases the cases considered in the present paper. 

Finally in section 5 for the Riemannian metrics on a 
surface we obtain the classification of geodesically 
equivalent pairs at a non-regular point (the point of 
bifurcation of the eigenvalues of the transition operator). 
Note that for generic pair of Riemannian metrics on a 
surface the set of points of their proportionality consists 
of isolated points. Therefore it is natural to consider the 
case when two Riemannian metrics on a surface are 
proportional in an isolated point. Some results of the 
global topological nature (namely, about the number of the 
points of proportionality for a pair of globally 
geodesically equivalent Riemannian metrics on a sphere) 
were obtained in \cite{mat}, but the local classification 
surprisingly was not done before. The canonical conformal 
structure on a surface, associated with a Riemannian 
metric, plays the crucial role in this classification. 
Using this conformal structure and Dini's Theorem 
(Levi-Civita's Theorem in the case of a surface), one can 
associate any pair of geodesically equivalent metrics on a 
surface, which are proportional in an isolated point $q_0$, 
with some (multiple-valued) analytic function in a 
neighborhood of $q_0$ with a singularity at $q_0$. The 
analysis of this singularity gives us the required 
classification (Theorems \ref{bift} and \ref {gendini}).

\section{Geodesic equivalence and orbital diffeomorphism of 
the extremal flows}
\setcounter{equation}{0} \indent

{\bf 2.1 Existence of the orbital diffeomorphism} Let $G_1$
and $G_2$ be two sub-Riemannian metrics on a distribution 
$D$ of a manifold $M$, ${||\cdot||_1}_q$ and 
${||\cdot||_2}_q$ be the corresponding Euclidean norms on 
$D(q)$, $h_1$ and $h_2$ be the Hamiltonians, defined by 
(\ref{hamsR}), where $||\cdot||_q$ is replaced by 
${||\cdot||_i}_q$, $i=1,2$. Also let $H_1$ and $H_2$ be the 
$\frac{1}{2}$-level sets of $h_1$ and $h_2$ respectively, 
i.e. 
\begin{equation}
\label{level} H_i=\{\lambda\in T^*M: 
h_i(\lambda)=\frac{1}{2}\} . 
\end{equation}
Besides, for given distribution $D$ and metric $G$ on it
denote by $J^k(D,G)$ the space of $k$-jets of all $C^k$
curves admissible to $D$ and parametrized by length w.r.t.
the metrics $G$. By definition the $1$-jet $J^1(D,G)$
satisfies $$J^1(D,G)=\{(q,v):q\in M, v\in D(q),
||v||_q=1\}.$$ For given curve $\gamma$ we will denote by
$j_{t_0}^{(k)}\gamma$ the $k$-jet of the curve $\gamma$ at
the point $t_0$.
\begin{prop}
\label{sufdif} If for some neighborhood $U$ of $q_0$ in $M$
there exist a fiberwise diffeomorphism $\Phi:H_1\cap
T^*U\mapsto H_2\cap T^*U$ and a function $a:H_1\mapsto 
\mathbb{R}$ such that 
\begin{equation} \label{conjrel} \Phi_*\vec h_1(\lambda)=
a(\lambda) \vec h_2\bigl(\Phi(\lambda)\bigr),
\end{equation} then the metrics $G_1$ and $G_2$ are
geodesically equivalent at $q_0$.
\end{prop}

{\bf Proof.} Indeed, $\Phi$ maps any trajectory of the
system $\dot\lambda=\vec h_1(\lambda)$, lying in $H_1\cap
T^*U$, to the curve, which coincides, up to
reparametrization, with a trajectory of the system
$\dot\lambda=\vec h_2(\lambda)$. Therefore in $U$ any 
normal sub-Riemannian geodesics of $G_1$ is, up to 
reparametrization, a normal sub-Riemannian geodesics of 
$G_2$. $\Box$

In the case of Riemannian metrics the relation
(\ref{conjrel}) is also necessary for the geodesic 
equivalence of metrics $G_1$ and $G_2$. Indeed, in this 
case there is only one geodesic passing through the given 
point in the given direction, and the map 
$P_i^{(1)}:H_i\mapsto J^1(TM,G_i)$ defined by 
\begin{equation}
\label{jet1} P_i(\lambda)\stackrel{def}{=}
j_0^1\pi(e^{t\vec h_i})=\Bigl(\pi(\lambda),\pi_*\bigl(\vec
h_i(\lambda)\bigr)\Bigr),\quad i=1,2,
\end{equation}
 is a diffeomorphism (here we denote by  $\pi:T^*M\mapsto M$ the canonical projection and by $e^{t \vec h_i}$ the flow generated by vector
field $\vec h_i$,
 $i=1,2$). So, directly by definition,
if the metrics $G_1$ and $G_2$ are geodesically equivalent
at $q_0$, then there is a neighborhood $U$ of $q_0$ such
that the following diffeomorphism
\begin{equation}
\label{orbriem} \Phi(\lambda)=
\bigl(P_2^{1}\bigr)^{-1}\left(\frac{1}{{||P_1^{1}(\lambda)||_2}_q}
P_1^{1}(\lambda)\right), \quad q=\pi(\lambda),
\end{equation}
 is fiberwise, maps $H_1\cap
T^*U$ to $H_2\cap T^*U$ and
 satisfies
(\ref{conjrel}) on $H_1\cap T^*U$.
\begin{defin}
A fiberwise diffeomorphism $\Phi$ defined on a nonempty
open set ${\mathcal B}$ of $H_1$ such that $\Phi({\mathcal
B})\subset H_2$ is called the orbital diffeomorphism  
of the extremal flows of the sub-Riemannian metrics $G_1$
and $G_2$ on ${\mathcal B}$, if it satisfies
(\ref{conjrel}) for any $\lambda\in {\mathcal B}$.
\end{defin}

Let us study the question, whether the existence of the 
orbital diffeomorphism is necessary for the geodesic 
equivalence of sub-Riemannian metrics. In the case of a 
proper sub-Riemannian metric (i.e., $D\neq TM$, $D$ is 
bracket-generating) an entire family of normal 
sub-Riemannian geodesics passes in general through the 
given point in the given direction. So, in order to 
distinguish different normal geodesics, passing through a 
point, we need jets of higher order. Besides, the presence 
of the abnormal extremal trajectories causes to addition 
difficulties, as shown below. By analogy with (\ref{jet1}) 
let us define the following mapping $P_i^{(k)}:H_i\mapsto 
J^k(D,G_i)$, $i=1,2$: 
\begin{equation}
\label{jetk}
P_i^{(k)}(\lambda)\stackrel{def}{=}j_0^k\pi(e^{t\vec h})
\end{equation}
Then one can check without difficulties that:

a) if $D$ is contact then the mapping $P_i^{(2)}$
establishes the diffeomorphism between $H_i$ and its image;

b) if $D$ is quasi-contact, ${\mathcal C}$ is the abnormal
line distribution of $D$, and the set ${\mathcal
S_i}\subset H_i$ is defined by
\begin{equation}
\label{sing} {\mathcal S}_i=\{\lambda\in
H_i:P^{1}_i(\lambda)\in {\mathcal C}\}, \end{equation} then
the restriction of the mapping $P_i^{(2)}$ on
${H_i}_q\backslash {\mathcal S}_i$ establishes the
diffeomorphism between ${H_i}_q\backslash {\mathcal S}_i$
and its image, while the restriction of $P_i^{(2)}$ on 
${\mathcal S}_i$ is constant on each fiber. 


Now denote by $\Omega_q (D, G_i)$ the set of all $C^\infty$
admissible curves, starting at $q$ and parametrized by
length w.r.t. the metric $G_i$ and let $J^k_q(D, G_i)$ be
the space of $k$-jet of these curves at $0$. Consider the
mapping $I_q: \Omega_q (D, G_1)\mapsto \Omega_q (D, G_2)$
which sends a curve $\gamma$ to its reparametrization
(w.r.t. the length of $G_2$). Obviously , this mapping
induces the diffeomorphisms $I_q^{(k)}: J^k_q(D,
G_1)\mapsto J^k_q(D, G_2)$. Collecting all such
diffeomorphisms for any $q$ we obtain a diffeomorphism
$I^{(k)}: J^k(D, G_1)\mapsto J^k(D, G_2)$. Then similarly
to (\ref{orbriem}) we obtain that if the distribution $D$
is one of the two listed in Introduction, and the
sub-Riemannian metrics $G_1$ and $G_2$, defined on $D$, are
geodesically equivalent at $q_0$, then there exist a
neighborhood $U$ of $q_0$ such that the following mapping
\begin{equation}
\label{orbsr} \Phi(\lambda)= 
\bigl(P_2^{(2)}\bigr)^{-1}\circ I^{(2)}\circ 
P_1^{(2)}(\lambda), 
\end{equation}
is well defined on the set ${\mathcal B}$, 
where ${\mathcal B}=H_1 \cap T^*U$ in contact case and
${\mathcal B}=(H_1 \cap T^*U)\backslash {\mathcal S}_1$ in
quasi-contact 
(here ${\mathcal S}_1$ is as in (\ref{sing})). Moreover, 
such $\Phi$ is the orbital diffeomorphism on the set 
${\mathcal B}$ w.r.t. the metrics $G_1$ and $G_2$. We have 
proved the following 
\begin{prop}
\label{subrorb}If $G_1$ and $G_2$ are Riemannian metric or
sub-Riemannian metrics defined on contact or quasi-contact 
distributions and if they are geodesically equivalent at 
some point $q_0$, then for some neighborhood $U$ of $q_0$ 
there exists the orbital diffeomorphism of the extremal 
flows of the metrics $G_1$ and $G_2$ on some nonempty open 
set ${\mathcal B}$ in $H_1 \cap T^*U$, $\pi(\mathcal B)=U$. 
In the Riemannian and contact case one can take ${\mathcal 
B}=H_1 \cap T^*U$, while in quasi-contact 
case one can take ${\mathcal B}= 
(H_1 \cap T^*U)\backslash {\mathcal S}_1$, where ${\mathcal 
S}_1$ is as in (\ref{sing}). 
\end{prop}

Actually, there is an analogue of the previous proposition 
for sub-Riemannian metrics defined on much more wide class 
of distributions. To formulate it let us introduce some 
notations. Denote by ${\mathcal A}_{q_0}(D)$ the set of all 
points $q\in M$ which can be connected with $q_0$ by 
abnormal extremal trajectory of the distribution $D$. 
For example, in Riemannian and contact case $\mathcal
A_{q_0}(D)$ is empty; in quasi-contact 
case ${\mathcal A}_{q_0}(D)$ is the set $L_{q_0}\backslash 
\{q_0\}$, where $L_{q_0}$ is the leaf of the abnormal line 
distribution, passing through $q_0$. 

\begin{prop}
\label{abnprop} Suppose that the sub-Riemannian metrics 
$G_1$ and $G_2$, defined on the bracket-generating 
distribution $D$, are geodesically equivalent at the point 
$q_0$ and for any neighborhood $V$ of $q_0$ the set 
$V\backslash{\mathcal A}_{q_0}(D)$ has positive Lebesgue 
measure. Then for some neighborhood $U$ of $q_0$ there 
exists the orbital diffeomorphism of the extremal flows of 
the metrics $G_1$ and $G_2$ on some open set ${\mathcal B}$ 
in $H_1 \cap T^*U$, $\pi(\mathcal B)=U$. 
\end{prop}

\begin{remark}
\label{minrem}
 Actually in the previous proposition one can replace  the set ${\mathcal
A}_{q_0}(D)$ by the set of all points $q\in M$ which can be
connected by abnormal extremal trajectory, having minimal
length w.r.t. the metric $G_1$ (or $G_2$) among all
admissible trajectories with endpoints $q_0$ and $q$.
\end{remark}

Since in the present paper we solve completely the problem 
of geodesic equivalence only in the cases, considered in  
Proposition \ref{subrorb}, we postpone the proof of 
Proposition \ref{abnprop} and the statement in Remark 
\ref{minrem} to the future paper. 
\medskip

{\bf 2.2 The orbital diffeomorphism in terms of the adapted 
frame to the pair of metrics.} Suppose that $D$ is an 
$(m,n)$-distribution on a manifold $M$. Let $q_0$ be a 
regular point w.r.t. the metric $G_1$ and $G_2$ (see 
Definition \ref {regdef}). It is simple to show that the 
regularity of the point $q_0$ is equivalent to the fact 
that the set of the multiplicities of the eigenvalues of 
the transition operator $S_q$ is the same for all points 
$q$ from some neighborhood of $q_0$. Therefore in some 
neighborhood $U$ of $q_0$ one can choose the basis 
$(X_1,\ldots , X_m)$ of the distribution $D$ orthonormal 
w.r.t. the metric $G_1$ such that each $X_i(q)$ is 
eigenvector of the transition operator $S_q$, $q\in V$. 
Such basis of $D$ will be called the {\it adapted basis to 
the ordered pair of metrics $(G_1,G_2)$ on a set $U$}. A 
frame $(X_1,\ldots, X_n)$ will be called the {\it adapted 
frame to the ordered pair of sub-Riemannian metrics 
$(G_1,G_2)$ on a set $U$}, if the tuple $\bigl(X_1,\ldots, 
X_m)$ is the adapted basis of $D$ w.r.t. $(G_1,G_2)$ on 
$U$.

Let us express the relation (\ref{conjrel}) for the orbital
diffeomorphism in terms of some adapted frame $(X_1,\ldots,
X_n)$ . Let $u_i:T^*M\mapsto\mathbb R$ be the
"quasi-impulse" of the vector field $X_i$,
\begin{equation}
\label{quasi} u_i(p,q)=p\bigl(X_i(q)\bigr),\quad q\in U,
p\in T^*U. \end{equation} For given diffeomorphism $\Phi$
defined on an open set of $T^*M$ denote by \begin{equation}
\label{comp} \Phi_i=u_i\circ \Phi,\quad 1\leq i\leq n.
\end{equation}
Suppose also that for any $i$, $1\leq i\leq m$, the 
eigenvalue of the transition operator $S_q$, corresponding 
to the eigenvector $X_i(q)$, is equal to $\alpha_i^2(q)$. 

\begin{lemma}
\label{firstcomp} If $\Phi$ is the orbital diffeomorphism
 of the extremal flows of  the metrics
$G_1$ and $G_2$ on an open set $\mathcal{B}\subset H_1\cap 
T^*U$, then the functions $\Phi_i$ with $1\leq i\leq m$ 
satisfy 
\begin{equation} \label{1compeq} \Phi_i=\frac{\alpha_i^2
u_i }{\sqrt{\sum_{k=1}^m\alpha_k^2 u_k^2}},\quad 1\leq
i\leq m.
\end{equation}
\end{lemma} {\bf Proof.}
 Since by construction the tuple
$(X_1,\ldots, X_m)$ constitute an orthonormal basis of the
distribution $D$ w.r.t. the metric $G_1$, the Hamiltonian
$h_1$ satisfies $h_1=\frac{1}{2}\sum_{i=1}^mu_i^2$, and
\begin{equation}
\label{Ham} \vec h_1=\sum_{i=1}^m u_i\vec u_i, \quad 
\pi_*\vec h_1=\sum_{i=1}^m u_iX_i ,\quad 
H_1=\Bigl\{\lambda\in T^*U:\sum_{i=1}^m u_i^2=1\Bigr\} 
\end{equation}
(here $\pi:T^*M\mapsto M$ is the canonical projection).
Let $\bar X_i$ be
\begin{equation}
\label{barXdef} \bar X_i=\frac{1}{\alpha_i}X_i,\quad 1\leq
i \leq m,
\end{equation}
and $\bar u_i(p,q)=p\bigl(\bar X_i(q)\bigr)$ be the corresponding
quasi-impulses. Then
\begin{equation}
\label{baru} \bar u_i=\frac{u_i}{\alpha_i},\quad 1\leq
i \leq m,
\end{equation}
Note that by construction $(\bar X_1,\ldots, \bar X_n)$ is the orthonormal basis of $D$ w.r.t. the metric $G_2$. Hence, similarly to (\ref {Ham}), we have
$\vec h_2=\sum_{i=1}^m \bar u_i\vec {\bar u}_i$, which together with
(\ref{barXdef}) and (\ref{baru}) implies that
\begin{equation}
\label{Ham2} \pi_*\vec h_2=\sum_{i=1}^m 
\frac{u_i}{\alpha_i^2}X_i , \quad H_2=\Bigl\{\lambda\in 
T^*U: \sum_{i=1}^m\bar u_i^2=1\Bigr\}= \Bigl\{\lambda\in 
T^*U: \sum_{i=1}^m \frac{u_i^2}{\alpha_i^2}=1\Bigr\} 
\end{equation}
Suppose that $\Phi$ is the orbital diffeomorphism on some 
set $\mathcal B$, satisfying (\ref{conjrel}) for some 
function $a$. Then by definition $\Phi(\lambda)\in H_2$ for 
any $\lambda\in {\mathcal B}$. This together with 
(\ref{comp}) and (\ref{Ham2}) implies that
\begin{equation}
\label{phiH2}
\sum_{i=1}^m  \frac{ \Phi_i^2}{\alpha_i^2}=1.
\end{equation}
Further from the fact that  $\Phi$ is fiberwise and (\ref{Ham}) it follows that 
$$(\pi_*\circ \Phi_*) \vec h_1(\lambda)=\pi_*\vec h_1(\lambda)=\sum_{i=1}^m u_i X_i.$$
On the other hand, (\ref{comp}) and (\ref{Ham2}) imply
$$\pi_*\vec h_2\bigl(\Phi(\lambda)\bigr)=
\sum_{i=1}^m \frac{\Phi_i}{\alpha_i^2}X_i.$$
From the last two relations and (\ref{conjrel}) it follows that
$$a \Phi_i=\,\alpha_i^2 u_i,\quad 1\leq i\leq m$$
From this and (\ref{phiH2}) it follows easily that
\begin{equation}
\label{a}
a=\sqrt{\sum_{k=1}^m\alpha_k^2 u_k^2},
\end{equation}
which implies (\ref{1compeq}).
$\Box$
\medskip

Now we will find the relation for the remaining components 
$\Phi_i$, $m+1\leq i\leq n$, of $\Phi$. Let $c_{j i}^k$ be 
the structural functions of the adapted frame $(X_1,\ldots, 
X_n)$, i.e., the function, satisfying $[X_i, X_j]=\sum 
c_{ji}^kX_k$. Let the vector fields $X_i$, $1\leq i\leq m$, 
satisfy (\ref{barXdef}) and set 
\begin{equation}
\label{barXeq} \bar X_i=X_i, \quad m+1\leq i\leq n.
\end{equation}

Note that by construction $(\bar X_1,\ldots, \bar X_n)$ is 
the adapted frame w.r.t. the ordered pair $(G_2, G_1)$. Let 
$\bar c_{j i}^k$ be the structural functions of the frame 
$(\bar X_1,\ldots, \bar X_n)$. The following functions will 
be very useful in the sequel together with function $a$, 
defined by (\ref{a}): 
\begin{equation}
\label{R} R_j\stackrel{def}{=} \frac{1}{2}\vec
h_1(\alpha_j^2)u_j+\alpha_j^2\vec
h_1(u_j)-\frac{1}{2}\alpha_j^2u_j \frac{\vec
h_1(a^2)}{a^2}-\sum_{1\leq i,k\leq m}\bar
c_{ji}^k\alpha_i\alpha_j\alpha_k u_iu_k,
\end{equation}
\begin{equation}
\label{Q} Q_{jk} \stackrel{def}{=} \sum_{i=1}^m \bar
c_{ji}^k\alpha_i u_i
\end{equation}

\begin{lemma}
\label{lastcomp} A map $\Phi$ is the orbital diffeomorphism
on a set $\mathcal{B}$ of the extremal flows of the metrics
$G_1$ and $G_2$ iff on $\mathcal B$ the functions $\Phi_k$
with $m+1\leq k\leq n$ satisfy the following relations:
\begin{equation}
\label{I} \forall 1\leq j\leq m:\quad \alpha_j
\sum_{k=m+1}^n
Q_{jk}
\Phi_k=\frac{R_j}{a}
\end{equation}
\begin{equation}
\label{II}
 \forall m+1\leq s\leq n:\quad
 \vec h_1(\Phi_s)- \sum_{k=m+1}^n
Q_{sk}\Phi_k =
\frac{1}{a}
\sum_{k=1}^m Q_{sk}\alpha_k u_k.
\end{equation}
\end{lemma}

{\bf Proof.}
In the sequel we set
\begin{equation}
\label{alphag}
\forall m+1\leq n:\quad  \alpha_i\equiv 1.
\end{equation}
Denote by $Y_i$ the vector field on $H_1$, which is the
lift of the vector field $X_i$ (i.e., $\pi_* Y_i=X_i$), and
$du_j(Y_i)=0$ $\forall 1\leq j\leq n$ (i.e., $Y_j$ is
horizontal field of the connection on $T^*M$ defined by 
distribution, satisfying $d u_1=\ldots = d u_n=0$). 
Similarly, let $\bar Y_i$ be the vector field on $H_2$, 
which is the lift of $\bar X_i$ and $d\bar u_j(Y_i)=0$ for 
all $1\leq j\leq n$. Note also that the tuple $(u_1,\ldots, 
u_n)$ defines the coordinates on each fiber $T^*_q M$ of 
$T^*M$. So, one can define the vector fields $\partial_ 
{u_i}$, $1\leq i\leq n$, as follows: $\partial_{u_i}$ is 
vertical (i.e., tangent to the fibers of $T^*M$) and 
$du_j(\partial {u_i})=\delta_{ij}$ for all $j=1,\ldots n$, 
where $\delta_{ij}$ is the Kronecker symbol. In the same 
way one can define the fields $\partial_{\bar u_i}$. With 
this notations, using (\ref{barXdef}) and (\ref{baru}), 
$\forall 1\leq i\leq n$ one can easily obtain the following 
relation 
: 
\begin{equation}
\label{barY}
\partial_{\bar u_i}=\alpha_i \partial_{u_i},\quad
\bar Y_i=\frac{1}{\alpha_i}
\left(Y_i+\sum_{j=1}^m\frac {X_j(\alpha_i)}{\alpha_i}u_j\partial_{u_j}\right)
\end{equation}
Besides, by standard calculations, we have
\begin{equation}
\label{h1str}
\vec h_1=\sum_{i=1}^m u_i\vec u_i=\sum_{i=1}^mu_i
Y_i+\sum_{i=1}^m\sum_{j,k=1}^nc_{ji}^ku_i u_k
\partial_{u_j},
\end{equation}
\begin{equation}
\label{h2str}
\vec h_2=\sum_{i=1}^m \bar u_i\vec {\bar u}_i=\sum_{i=1}^m\bar u_i
\bar Y_i+\sum_{i=1}^m\sum_{j,k=n}^m\bar c_{ji}^k\bar u_i \bar u_k
\partial_{\bar u_j}.
\end{equation}

Substituting (\ref{baru}) and (\ref{barY}) into (\ref{h2str}), we obtain
\begin{equation}
\label{h2sub}
\vec h_2=\sum_{i=1}^m \frac{u_i}{\alpha_i^2} Y_i+\sum_{i,j=1}^m \frac{X_i(\alpha_j)}{\alpha_i^2\alpha_j}u_i u_j\partial_{u_j}+
\sum_{i=1}^m\sum_{j,k=1}^n\frac{\bar c_{ji}^k\alpha_j}{\alpha_i\alpha_k}
u_i u_k
\partial_{u_j}.
\end{equation}
This together with (\ref{1compeq}) implies easily that
\begin{eqnarray}
 &~& \vec
h_2\bigl(\Phi(\lambda)\bigr)=a^{-1}\, \sum_{i=1}^mu_i 
Y_i+a^{-2}\,\sum_{j=1}^m \Bigl(\frac{1}{2}\vec 
h_1(\alpha_j^2)u_j+\sum_{i,k=1}^m \bar 
c_{ji}^k\alpha_i\alpha_j\alpha_k u_i 
u_k\Bigr)\partial_{u_j}+\nonumber\\ &~& a^{-1} \sum_{j=1}^m 
\sum_{k=m+1}^n\sum_{i=1}^m \bar c_{ji}^k\alpha_i\alpha_j 
u_i \Phi_k \partial_{u_j}+\sum_{s=m+1}^n 
\sum_{i=1}^m\Bigl(a^{-2}\sum_{k=1}^m\bar 
c_{si}^k\alpha_i\alpha_k u_i u_k+\Bigr.\label{h2fin}\\ 
&~&\Bigl. a^{-1}\sum_{k=m+1}^n \bar c_{si}^k\alpha_i u_i 
\Phi_k\Bigr)\partial_{u_s}, 
\end{eqnarray}
where $a$ is as in (\ref{a}).
On the other hand, from the fact that $\Phi$ is fiberwise and relation (\ref{1compeq}) it follows that
\begin{equation}
\label{h1fin}
\Phi_*\vec h_1(\lambda)=\sum_{i=1}^m u_i Y_i+
\sum_{j=1}^m\vec h\left(\frac{\alpha_j^2
u_j }{\sqrt{\sum_{l=1}^m\alpha_l^2 u_l^2}}\right)\partial_{u_j}+
\sum_{j=m+1}^n\vec h(\Phi_j)\partial_{u_j}.
\end{equation}
Using relations (\ref{h2fin}) and (\ref{h1fin}), it is not
hard to check by direct calculations that (\ref{conjrel})
holds iff both (\ref{I}) and (\ref{II}) hold, which
concludes the proof of the Lemma. $\Box$
\medskip


{\bf 2.3 The first divisibility condition.} Let ${\mathcal
I}_1:D(q)^*\mapsto D(q)$ be the canonical isomorphism
w.r.t. the inner product ${G_1}_q(\cdot,\cdot)$, namely,
$\ell(\cdot)={G_1}_q({\mathcal I}_1(\ell),\cdot)$ $\forall
\ell\in D(q)^*$. Define the following function ${\mathcal
P}:T^*M\mapsto \mathbb R$:
\begin{equation}
\label{funP} {\mathcal P}(p,q)=\Bigl(||{\mathcal
I}_1(p|_{D(q)})||_{_{2 _q}}\Bigr)^2, \quad q\in M, p\in
T_q^*M
\end{equation}
(here $||\cdot||_{2_q}$ is the Euclidean norm on $D(q)$
corresponding to the inner product ${G_2}_q(\cdot,\cdot)$,
$p|_{D(q)}$ is the restriction of covector $p\in T_q^*M$ on
the subspace $D(q)$). Obviously, the restriction of
${\mathcal P}$ on each fiber $T^*_qM$ is a degree 2
homogeneous polynomial, while the restriction of $\vec h_1(
{\mathcal P})$ on each fiber $T^*_qM$ is a degree 3
polynomial. Besides, in a neighborhood of the regular point 
\begin{equation}
\label{Pa} {\mathcal P}=a^2=\sum_{i=1}^m\alpha_i^2 u_i^2,
\end{equation}
where $(u_1,\ldots, u_m)$ are quasi-impulses of the vectors
of the adapted basis $(X_1,\ldots X_m)$ to the order pair
$(G_1,G_2)$ and $\alpha_i^2$ are eigenvalues of the
transition operator $S_q$, corresponding to the
eigenvectors $X_i$.
\begin{defin}
\label{1divcond} We will say that the ordered pair
$(G_1,G_2)$ of sub-Riemannian metrics on the distribution
$D$ satisfies the first divisibility condition on a set
$U$, if the polynomial $\vec h_1( {\mathcal P})|_{T_q^*M}$
is divided by the polynomial ${\mathcal P}|_{T_q^*M}$ for
any $q\in U$.
\end{defin}

\begin{prop}
\label{1divpr} Let $D$ be an $(m,n)$-distribution on a 
manifold $M$ such that 
\begin{equation}
\label{jump1} \forall 1\leq s\leq n-m+1,\quad \dim
D^s=m+s-1.
\end{equation}
Suppose also that for given two sub-Riemannian metrics 
$G_1$ and $G_2$ on $D$ and for some open set $U$ of $M$ 
there exists an orbital diffeomorphism of the extremal 
flows of these metrics in some open set $\mathcal B$ in 
$H_1\cap T^*U$, $\pi(\mathcal B)=U$. 
Then the pair $(G_1, G_2)$ satisfies the first divisibility
condition on $U$.
\end{prop}

{\bf Proof.}
Since the set of regular points is dense ( Proposition 
\ref{denseprop}) it is sufficient to prove the first 
divisibility condition for a regular point $q_0$ w.r.t. the 
pair $(G_1,G_2)$. Therefore in order to obtain the first 
divisibility condition we can use Lemmas \ref{firstcomp} 
and \ref{lastcomp}. Note also that by (\ref{R}) the 
function $R_j$ has the following form on each fiber: 
\begin{equation}
\label{Rp} R_j=-\frac{1}{2}\alpha_j^2u_j\frac{\vec 
h_1({\mathcal P})}{{\mathcal P}}+{\rm polynomial} 
\end{equation}

First suppose that $D=TM$ (in this case the assumption
(\ref{jump1}) holds automatically). Then the identity
(\ref{I}) is equivalent to the identity $R_j\equiv0$, 
$1\leq j\leq n$, which holds on open set $\mathcal B$ in 
$H_1\cap T^*U$ with $\pi(\mathcal B)=U$ and therefore on 
the whole $T^*U$. 
Hence from (\ref{Rp})
it follows that
$u_j\frac{h_1({\mathcal P})}{{\mathcal P}},$ has to be a
polynomial, which implies easily that the polynomial
${\mathcal P}$ has to divide the polynomial $\vec
h_1({\mathcal P})$, i.e., the first divisibility condition
holds.

Now consider the case $D\neq TM$. By assumption
(\ref{jump1}), we can complete the adapted basis
$(X_1,\ldots, X_m)$ to the adapted frame such that
\begin{equation}
\label{compjump} \forall m+1\leq s\leq n \quad X_s\in
D^{s-m+1}
\end{equation}
Then $D^2={\rm span}(X_1,\ldots X_{m+1})$, which implies
that there exist indices $\bar i,\bar j$, $1\leq \bar i, 
\bar j\leq m$, such that 
 $\bar c_{\bar j\bar i}^{m+1}(q_0)\neq 0$, while $\bar c_{ij}^k=0$ for all
 $1\leq i,j\leq m$ and $k>m+1$. In other words,
 \begin{eqnarray}
&~&\label{Qj0} \forall k,j:\, k>m+1, 1\leq j\leq m \quad
Q_{jk}\equiv 0\\ &~& \label{Qjn0} \exists \bar j: 1\leq
\bar j\leq m , \quad Q_{\bar j m+1}\not\equiv 0
\end{eqnarray}
(see (\ref{Q}) for the definition of the functions
$Q_{jk}$). Then from (\ref{I}) it follows that
\begin{equation}
\label{Philast} \Phi_{m+1}=\frac{R_{\bar j}}{\alpha_{\bar
j}Q_{\bar {j} m+1}\sqrt{\mathcal P}}.
\end{equation}
Using (\ref{Rp}), we obtain
\begin{equation}
\label{Philast1} \Phi_{m+1}=-\frac{1}{2}\alpha_{\bar
j}u_{\bar j}\frac{\vec h_1({\mathcal P})}{Q_{{\bar
j}m+1}{\mathcal P}^{3/2}}+\frac{1}{\alpha_{\bar j}Q_{{\bar
j} m+1}\sqrt{\mathcal P}}\,\,{\rm polynomial}
\end{equation}
on each set $\mathcal B\cap T_q^*M$, $q\in U$.

Further, from assumption (\ref{jump1}) it follows that
 $\bar c_{si}^{k+1}=0$ for any $k,s,i$ such that $m<s<k$ and  $1\leq
i\leq m$. On the other hand, there exist $\bar i$, $1\leq
\bar i\leq m$, such that $\bar c_{s\bar i}^{s+1}\neq 0$. In
other words,
\begin{eqnarray}
&~&\label{Q0} \forall k,\, s:
m<s<k \quad Q_{s\,k+1}\equiv 0\\ &~& \label{Qn0} \forall s:
m<s<n-1\quad Q_{s\,s+1}\not\equiv 0.
\end{eqnarray}
Hence by (\ref{II}), applied for $s=m+l-1$ with $2\leq
l\leq n-m$, one has
\begin{equation}
\label{Phiprel} \Phi_{m+l}= Q_{m+l-1\, m+l}^{-1}\Bigl(\vec
h_1(\Phi_{m+l-1})- \sum_{k=m+1}^{m+l-1}
Q_{m+l-1\,k}\Phi_k -
{\mathcal P}^{-1/2}
\sum_{k=1}^m Q_{m+l-1,k}\alpha_k u_k\Bigr).
\end{equation}
Then by induction from (\ref{Philast1}) and (\ref{Phiprel})
it is not difficult to get the following relation for any
$2\leq l\leq m-n$
\begin{equation}
\label{Phiind} \Phi_{m+l}=\frac{(-1)^l(2l-1)!!\, 
\alpha_{\bar j} u_{\bar j}\bigr(\vec h_1({\mathcal 
P})\bigl)^l}{2^l Q_{\bar j\, m+1}{\mathcal 
P}^{l+1/2}\prod_{i=1}^{l-1} Q_{m+i\,m+i+1}}+\frac{{\rm 
polynomial}}{Q_{\bar j\,m+1}^l{\mathcal 
P}^{l-1/2}\prod_{i=1}^{l-1} Q_{m+i\,m+i+1}^{l-i}} 
\end{equation}
on each $\mathcal B\cap T_q^*M$, $q\in U$. Substituting the 
expression for $\Phi_{m+l}$ from (\ref{Phiind}) to identity
(\ref{II}) with $s=n$ one can obtain without difficulties 
that \begin{equation} \label{finid} \cfrac{u_{\bar 
j}\bigr(\vec h_1({\mathcal P})\bigl)^{n-m+1}}{ Q_{\bar 
j\,m+1}{\mathcal P}^{n-m+3/2} \prod_{i=1}^{n-m-1} 
Q_{m+i\,m+i+1}}=\cfrac{{\rm polynomial}}{Q_{\bar 
j\,m+1}^{n-m+1}{\mathcal P}^{n-m+1/2}\prod_{i=1}^{n-m-1} 
Q_{m+i\,m+i+1}^{n-m-i+1}} 
\end{equation}
or, equivalently
 \begin{equation} \label{finid1} \cfrac{u_{\bar
j}Q_{\bar j\,m+1}^{n-m}\bigl(\prod_{i=1}^{n-m-1}
Q_{m+i\,m+i+1}^{n-m-i}\bigr)\bigl(\vec h_1({\mathcal
P})\bigr)^{n-m+1}}{\mathcal P}={\rm polynomial}.
 \end{equation}
 on each set $\mathcal B\cap T_q^*M$, $q\in U$. Note that
 the left-hand side of (\ref{finid1}) is rational
function. Hence from (\ref{finid1}) it has to be
polynomial. Note also that ${\mathcal P}$ is positive 
definite quadratic form (see (\ref{Pa})), while the 
functions $Q_{j k}$ are linear (with real coefficients) on 
each fiber. Therefore from (\ref{finid1}) it follows easily 
that the polynomial ${\mathcal P}$ has to divide the 
polynomial $\vec h_1({\mathcal P})$, i.e., the first 
divisibility condition holds. The proof of the proposition 
is concluded. $\Box$. 
\medskip

Note that if $D$ is contact, quasi-contact,
or
$D=TM$, then the assumption (\ref{jump1}) of the previous
proposition holds. So, as a direct consequence of
Proposition \ref{subrorb} and the previous proposition, we
have the following
\begin{cor}
\label{divcor1}
 Suppose that two metrics $G_1$ and $G_2$
defined on the distribution $D$ are geodesically equivalent
at the point $q_0$. Assume also that the distribution $D$
satisfies one of the two following conditions:
\begin{enumerate} \item $D=TM$ (the Riemannian case); \item
$D$ is corank 1 contact or quasi-contact distribution;
\end{enumerate}
Then the pair $(G_1, G_2)$ satisfies the first divisibility
condition on $U$.
\end{cor}

So, in the cases under consideration the first divisibility
condition is necessary for the geodesic equivalence. In the
next proposition we collect all information from the first
divisibility condition, which will be used in the sequel.
It shows that the first divisibility condition imposes
rather strong restrictions on the pair of the metrics.
\begin{prop}
\label{divcor} Suppose that the metrics $G_1$ and $G_2$,
defined on the distribution $D$, satisfy the first
divisibility condition on some set $U$. If $(X_1,\ldots
X_m)$ is a basis of $D$ adapted to the order pair
$(G_1,G_2)$, and the transition operator $S_q$ has the form
$S_q={\rm diag}\,(\alpha_1^2(q),\ldots, \alpha_m^2(q))$ in
this basis ($\alpha_i>0$), then the following relations
hold
\begin{eqnarray}
&~& \label{cor4} [X_i,X_j](q)\notin D(q)\,\,\Rightarrow\,\, 
\alpha_i(q)=\alpha_j(q); \\ &~&\label{cor1} 
X_i\left(\frac{\alpha_j^2}{\alpha_i^2}\right)=2c_{ji}^j 
\left(1-\frac{\alpha_j^2}{\alpha_i^2}\right); 
\\
&~&\label{cor2}
X_i\left(\frac{\alpha_j^2}{\alpha_i}\right)=0,\quad
\alpha_i\neq \alpha_j\\
&~&\label{cor3} 
X_i\left(\frac{\alpha_j}{\alpha_k}\right)=0,\quad 
\alpha_j\neq \alpha_i, \alpha_k\neq \alpha_i;\\ 
&~&\label{cor5} (\alpha_j^2-\alpha_i^2)c_{ji}^k+ 
(\alpha_j^2-\alpha_k^2)c_{jk}^i+(\alpha_i^2-\alpha_k^2)c_{ik}^j=0, 
\quad i,j,k\,\,{\rm are}\,\,{\rm pairwise}\,\,{\rm 
distinct}. 
\end{eqnarray}
(in all relations above $1\leq i,j,k \leq m$).
\end{prop}

{\bf Proof.} As before let us complete the adapted basis 
$(X_1,\ldots X_m)$ of $D$ somehow to the local frame. From 
(\ref{h1str}) and (\ref{Pa}) by direct calculation one has 
\begin{equation} \label{rel1} \vec
h_1({\mathcal P})=\sum_{i,j=1}^mX_i(\alpha_j^2) u_i
u_j^2+2\sum_{i,j=1}^m\sum_{k=1}^n c_{ji}^k\alpha_j^2 u_i
u_j u_k
\end{equation}
On the other hand by the first divisibility condition there 
exist functions $p_i(q)$, $1\leq i\leq n$, such that 
\begin{equation}
\label{rel2} \vec h_1({\mathcal P})=\Bigl(\sum_{i=1}^np_i
u_i\Bigr)\Bigl(\sum_{j=1}^m\alpha_j^2 u_j^2\Bigr).
\end{equation}

Relation (\ref{cor5}) follows immediately from comparing
the coefficients of $u_iu_ju_k$ in the right-hand sides of
(\ref{rel1}) and (\ref{rel2}), where $i,j,k$ are pairwise
distinct and $1\leq i,j,k\leq m$.

Further, comparing the coefficient of $u_iu_ju_k$ in the
right-hand side of (\ref{rel1})and (\ref{rel2}), where
$1\leq i\leq j\leq m$ and $k>m$, we have
\begin{equation}
\label{cor4d} c_{ji}^k(\alpha_i^2-\alpha_j^2)=0
\end{equation}
Therefore, if $[X_i,X_j](q)\notin D(q)$, then there exists
$k>m$ such that $c_{ji}^k(q)\neq 0$, which implies that
$\alpha_i(q)=\alpha_j(q)$. Relation (\ref{cor4}) is proved.

Further, comparing coefficients of $u_i^3$ in the
right-hand sides of (\ref{rel1}) and (\ref{rel2}) we obtain
that
\begin{equation}
\label{pi} p_i=\frac {X_i(\alpha_i^2)}{\alpha_i^2},
\end{equation}
while comparing coefficients of $u_i u_j^2$ with $i\neq j$
and using (\ref{pi}) one obtains easily that
\begin{equation}
\label{cor11} \alpha_i^2 X_i(\alpha_j^2)-\alpha_j^2
X_i(\alpha_i^2)=2c_{ji}^j(\alpha_i^2-\alpha_j^2)\alpha_i^2.
\end{equation}
The last equation implies (\ref{cor1}).

In order to prove (\ref{cor2}) note that we can obtain one
more relation in addition to (\ref{cor1}),
starting with the metric $G_2$ as the original one and
using transition from the metric $G_2$ to the metric $G_1$.
Namely, if
$\bar X_i$ is as in (\ref{barXdef})
 then by
analogy with (\ref{cor11}) we have
\begin{equation}
\label{barcor}\bar \alpha_i^2 \bar
X_i(\bar\alpha_j^2)-\bar\alpha_j^2 \bar
X_i(\bar\alpha_i^2)=2\bar
c_{ji}^j(\bar\alpha_i^2-\bar\alpha_j^2)\bar\alpha_i^2.
\end{equation}
Obviously, $\bar \alpha_i=\frac{1}{\alpha_i}$. Also, by
(\ref{barXdef}) one get easily that
\begin{equation}
\label{barcjij} \bar c_{ji}^j=\frac{c_{ji}^j}{\alpha_i}-
\frac{X_i(\alpha_j)}{\alpha_i\alpha_j}.
\end{equation}
 Substituting the
last two relations and (\ref{barXdef}) into (\ref{barcor})
it is not difficult to get the following
\begin{equation}
\label{invrel}
 \alpha_i^2 X_i(\alpha_j^2)-\alpha_j^2
X_i(\alpha_i^2)=2\alpha_j\bigl((\alpha_jc_{ji}^j-X_i(\alpha_j)\bigr)
(\alpha_i^2-\alpha_j^2)
\end{equation}
Then combining (\ref{cor11}) with (\ref{invrel}) and using
the fact that $\alpha_i\neq \alpha_j$ one get
\begin{equation}
\label{cjii}
c_{ji}^j=\frac{1}{2}\frac{X_i(\alpha_j^2)}{\alpha_j^2-\alpha_i^2}
\end{equation}
Substituting the last relation again in (\ref{cor11}) we
have
$$2\alpha_i^2X_i(\alpha_j^2)-\alpha_j^2X_i(\alpha_i^2)=0,$$
which is equivalent to (\ref{cor2}). Relation (\ref{cor3})
follows immediately from (\ref{cor2}).
\medskip

\begin{cor}
\label{rank2prop} If $D$ is a bracket-generating
$(2,n)$-distribution, $n>2$, and the metrics $G_1$ and
$G_2$, defined on $D$, satisfy the first divisibility
condition, then they are proportional, namely
$G_{2_{q}}=\alpha(q) G_{1_{q}}$. 
\end{cor}

{\bf Proof.} Since $D$ is bracket-generating, the set $V_1$
of points $q$ with
\begin{equation} \label{D23}
\dim D^2(q)=3
\end{equation}
is open and dense. By Proposition \ref{denseprop} the
intersection $V_2$ of this set with the set of all regular
points w.r.t. the metrics $G_1$ and $G_2$ is also open
dense. Therefore it is sufficient to proof the corollary
for the points of $V_2$. From regularity it follows the
existence of the adapted frame $(X_1, X_2)$. So, we can
apply the previous proposition. By (\ref{D23}),
$[X_1,X_2]\not\in D$.
Hence from (\ref{cor4}) it follows that
$\alpha_1\equiv\alpha_2$,
which completes the proof of the corollary.
 $\Box$

Suppose that $D$ is a step $1$ bracket-generating $(2, 
n)$-distribution ($\dim D^{l+1}=\dim D^{l}+1$ for any 
$1\leq l\leq n-m$). It can be shown easily that $D$ satisfy 
the assumptions of Proposition \ref{abnprop}. Therefore by  
Propositions \ref{abnprop}, \ref{1divpr}, and Corollary 
\ref{rank2prop} we have the following

\begin{prop}
\label{step1prop2} Suppose that two sub-Riemannian metrics 
$G_1$ and $G_2$ are defined on a step $1$ 
bracket-generating $(2,n)$-distribution, where $n>2$. If 
they are geodesically equivalent at some point $q_0$, then 
they are proportional, namely $G_{2_{q}}=\alpha(q) 
G_{1_{q}}$ in some neighborhood of $q_0$. 
\end{prop} 

Our conjecture is that the factor $\alpha(q)$ in the 
previous propoposition has to be constant, but we can prove 
it still only in the case $n=3$ (see Corollary \ref{23cor} 
below). We finish this section with the following useful 
lemma 

\begin{lemma}
\label{Rlemma} Suppose that the metrics $G_1$ and $G_2$,
defined on an $(m,n)$-distribution $D$, satisfy the first
divisibility condition at some neighborhood $U$ of a
regular point $q_0$. If $(X_1,\ldots X_m)$ is a basis of
$D$ adapted to the order pair $(G_1,G_2)$, and the
transition operator $S_q$ has the form $S_q={\rm
diag}\,(\alpha_1^2(q),\ldots, \alpha_m^2(q))$ in this basis
($\alpha_i>0$), then the functions $R_j$, $1\leq j\leq m$,
defined by (\ref{R}), can be written in the following form
\begin{eqnarray} &~& R_j=
\sum_{i=1} ^m(1-\delta_{ji})\Bigl((\alpha_j^2-\alpha_i^2) 
c_{ji}^i-\frac{X_j(\alpha_i^2)}{2}\Bigr)u_i^2+ 
\sum_{i=1}^m(1-\delta_{ji})\frac{\alpha_i^2}{2\alpha_j^2} 
X_i \Bigl(\frac{\alpha_j^4}{\alpha_i^2}\Bigr) 
u_iu_j+\nonumber \\ &~&\label{Rjexp} \sum_{i=1}^m\sum_{k=1} 
^m(1-\delta_{ik})(\alpha_j^2-\alpha_k^2)c_{ji}^ku_iu_k+ 
\alpha_j^2\sum_{i=1}^m\sum_{k=m+1}^nc_{ji}^ku_iu_k. 
\end{eqnarray}
(here $\delta_{ij}$ is the Kronecker symbol).
\end{lemma}
The relation can be obtain without difficulties by
substitution of (\ref{pi}), (\ref{barcjij}) and the
following obvious identity
\begin{equation}
 \label{relc2}
 \bar
 c_{i,j}^k=c_{i,j}^k\frac{\alpha_k}{\alpha_i\alpha_j}\quad
 i,j,k \,\,{\rm are}\,\,{\rm pairwise}\,\,{\rm distinct}
 \end{equation}
 into (\ref{R}).

\section {The case of Riemannian metrics near regular point}
\setcounter{equation}{0} \indent

In the present section, using the technique developed
above, we give a new proof of classical Levi-Civita's
Theorem about the classification of all Riemannian
geodesically equivalent metrics in a neighborhood of the
regular points w.r.t. these metrics (see \cite{levi},
\cite{top}). This proof is
 rather elementary and transparent from the geometrical point of view.
 Some crucial ideas of this proof will be used in the next section for obtaining the
 corresponding classification for sub-Riemannian geodesically
 equivalent metrics on quasi-contact distributions.

Here we prefer the coordinate-free formulation of
Levi-Civita's Theorem, which in our opinion clarifies the
statement of it. But before let us introduce some notations
and prove some preparatory lemmas.

Let $G_1$ and $G_2$ be Riemannian metrics on an 
$n$-dimensional manifold $M$. Let $q_0$ be a regular point 
w.r.t. these metrics. Suppose that $(X_1,\ldots X_n)$ is a 
frame adapted to the order pair $(G_1,G_2)$ in some 
neighborhood of $q_0$, and the transition operator $S_q$ 
from the metric $G_1$ to the metric $G_2$ has the form 
$S_q={\rm diag}\,(\alpha_1^2(q),\ldots, \alpha_n^2(q))$ in 
this basis ($\alpha_i>0$). 

Let $R_j$, $1\leq j\leq n$, be as in (\ref{R}). 
Propositions \ref{sufdif}, \ref{subrorb} and Lemma 
\ref{lastcomp} imply the following 

\begin{lemma}
\label{prepLC} Two Riemannian metrics $G_1$ and $G_2$ are
geodesically equivalent at a regular point $q_0$ if and
only if
 there exist some neighborhood $U$ of $q_0$
such that the following identities hold on $T^*U$
\begin{equation}
\label{Rj0} \forall j: 1\leq j\leq n\quad R_j\equiv0.
\end{equation}
\end{lemma}

Further, let $\{\lambda_1,\ldots,\lambda_N\}$ be the set of
all distinct eigenvalues of $S_q$, $\lambda_s>0$ (from the
regularity the number of these eigenvalues is constant for
all $q$ from some neighborhood of $q_0$). Denote by
\begin{equation}
\label{Is} I_s=\{i: \alpha_i^2=\lambda_s\}\quad 1\leq s\leq
N
\end{equation} Denote also by $D_s$ the following rank
$|I_s|$-distribution
\begin{equation}
D_s ={\rm span} \{X_i\}_{i\in I_s},\quad 1\leq s\leq N
\end{equation}

\begin{lemma}
\label{invcons} If two Riemannian metrics $G_1$ and $G_2$
are geodesically equivalent at a regular point $q_0$, then
the distribution $D_s$ is integrable in a neighborhood of
$q_0$.
\end{lemma}

{\bf Proof.} By Lemma \ref{prepLC} identities (\ref{Rj0}) 
hold. Taking the coefficient of $u_i u_k$ from 
(\ref{Rjexp}), by (\ref{Rj0}) one has 
\begin{equation}
\label{levi1}
(\alpha_j^2-\alpha_k^2)c_{ji}^k+(\alpha_j^2-\alpha_i^2)c_{jk}^i=0,\quad
i,j,k\,\,{\rm are}\,\,{\rm pairwise}\,\,{\rm distinct}.
\end{equation}
If $\alpha_i=\alpha_j$ and $\alpha_i\neq \alpha_k$, then
from the last relation $(\alpha_j^2-\alpha_k^2)c_{ji}^k=0$,
which implies that $c_{ji}^k=0$. In other words, if $i,j\in 
I_s$, then $[X_i, X_j]\in D_s$. So, $D_s$ is an integrable 
distribution. $\Box$ 
\medskip

\begin{lemma}
\label{twored} 
If two Riemannian metrics $G_1$ and $G_2$ are geodesically
equivalent at a regular point $q_0$, then the distribution
$$D_{s,l}\stackrel{def}{=} {\rm span}
\bigl(D_{s},D_{l}\bigr)={\rm span}\Bigl\{ X_i\Bigr\}_{i\in
(I_{s}\cup I_{l})},$$
is integrable in a neighborhood of $q_0$ for all $s,l$ ,
$1\leq s\neq l\leq N$.
\end{lemma}

{\bf Proof.} By the previous lemma it is sufficient to
prove that for any three indices $i,j,k$ with pairwise
distinct $\alpha_i$, $\alpha_j$, and $\alpha_k$ we have
$c_{ji}^k\neq 0$. Making the corresponding permutation of
indices in (\ref{levi1}), we obtain one more relation
\begin{equation}
\label{levi2}
-(\alpha_j^2-\alpha_k^2)c_{ji}^k+(\alpha_i^2-\alpha_j^2)c_{ik}^j=0,\quad
i,j,k\,\,{\rm are}\,\,{\rm pairwise}\,\,{\rm distinct}.
\end{equation}
Combining (\ref{levi1}), (\ref{levi2}), and (\ref{cor5}),
we obtain the system of three linear equations w.r.t.
$c_{ji}^k$, $c_{jk}^i$, and $c_{ik}^j$ with the determinant
equal to
$2(\alpha_i^2-\alpha_j^2)(\alpha_i^2-\alpha_k^2)(\alpha_k^2-\alpha_j^2)$,
which implies that $c_{ji}^k=0$. $\Box$ \medskip

 From the previous lemma by standard arguments one has the following

\begin{cor}
\label{coord} If two Riemannian metrics $G_1$ and $G_2$ are
geodesically equivalent at a regular point $q_0$,
then in some neighborhood $U$ of the point $q_0$ there
exist coordinates $(x_1,\ldots x_n)$
such that
\begin{equation}
\label{netcom}\forall s: 1\leq s\leq p\quad
D_s=\{d\,x_i=0\}_{i\not\in I_s}.
\end{equation}
In other words, in this coordinates the leaves of the
integrable distribution $D_s$ are $|I_s|$-dimensional
linear subspaces, parallel to the coordinate $\{x_i\}_{i\in
I_s}$-subspace.
\end{cor}

For any $s$, $1\leq s\leq N$, denote by ${\mathcal F}_s$
the foliation of the integral manifolds of the distribution
$D_s$. Let ${\mathcal F}_s(q_0)$ be the leaf of ${\mathcal
F}_s$, passing through the point $q_0$. Also, let $U$ be 
the neighborhood of $q_0$ from Corollary \ref{coord}. Then 
for any $s$, $1\leq s\leq N$, one can define a special map 
$pr_s:U\mapsto {\mathcal F}_s(q_0)$ in the following way: 
the point $pr_s(q)$ is the point of intersection of 
${\mathcal F}_s(q_0)$ with the integral manifold of the 
distribution ${\rm span} \{D_l: 1\leq l\leq N, l\neq s\}$, 
passing through $q$. In the coordinates of Corollary 
\ref{coord} with $q_0=(0,\ldots,0)$ the map $ pr_s$ is the 
projection on the coordinate $\{x_i\}_{i\in I_s}$-subspace 
(which preserves all coordinates $x_i$, ${i\in I_s}$). Now 
we are ready to formulate Levi-Civita's Theorem: 


\begin{theor}
(Levi-Civita) \label{dinit} Two Riemannian metrics $G_1$
and $G_2$
are geodesically equivalent at a point $q_0$ 
if and only if for any $s$, $1\leq s\leq N$, on a manifold
${\mathcal F}_s(q_0)$ there exist a Riemannian metric $g_s$ 
and a positive function $\beta_s$, which is constant if 
$\dim\, {\mathcal F}_s>1$, such that $\beta_s(q_0)\neq 
\beta_l(q_0)$ for all $s\neq l$ and in some neighborhood of 
$q_0$ the metrics $G_1$ and $G_2$ have the following form 
\begin{eqnarray}
\label{met1}
&~&
G_1=\sum_{s=1}^N \gamma_s (pr _s)^* g_s,\quad
\\ &~&
G_2=\sum_{s=1}^N\lambda_s \gamma_s (pr _s)^* g_s
\label{met2},
\end{eqnarray}
where
\begin{eqnarray}
&~&\label{lambdas}
\lambda_s=(\beta_s\circ pr_s)\prod_{l=1}^N(\beta_l\circ pr_l),\\
&~&
\label{gammas} \gamma_s=\prod_{l\neq
s}\Bigl|\frac{1}{(\beta_l\circ pr_l)}-\frac{1}{(\beta_s\circ pr_s)}\Bigr|
.
\end{eqnarray}
\end{theor}

{\bf Proof.} We start with the proof of the "only if" part. 
Below we work in the coordinate neighborhood $U$ of 
Corollary \ref{coord}. First let us prove the following 
\begin{lemma}
\label{lemma1} For any $s$, $1\leq s\leq N$, there exist a 
metric $g_s$ on ${\mathcal F}_s$ and some function
$\gamma_s$ such that (\ref{met1}) holds. 
\end{lemma}

{\bf Proof.} Since by construction for any $s_1\neq s_2$ 
the distributions $D_{s_1}$ and $D_{s_2}$ are orthogonal 
w.r.t. the metric $G_1$, the relation (\ref{met1}) is 
equivalent to the fact that for any $s$, $1\leq s\leq N$, 
there exists the metric $g_s$ on ${\mathcal F}_s$ and the 
function $\gamma_s$ such that 
\begin{equation}
\label{metorth} \forall Y\in D_s(q)\quad
{G_1}_q=\gamma_s{g_s}_{pr_s q}\bigl((d (pr_s)_ q Y\bigr)
\end{equation}
If $\dim {\mathcal F}_s=1$ (or, equivalently, $|I_s|=1$), 
then relation (\ref{metorth}) holds automatically for some 
$g_s$ on ${\mathcal F}_s$ and some function $\gamma_s$, 
because all quadratic forms of one variable are 
proportional. Let us prove (\ref{metorth}) 
in the general case. First, as $g_s$ we can take the 
restriction $G_1\Bigl|_{{\mathcal F}_s(q_0)}\Bigr.$ of 
$G_1$ to $\mathcal F_s(q_0)$, i.e., 
\begin{equation}
\label{gs} g_s= G_1\Bigl|_{{\mathcal F}_s(q_0)}\Bigr.
\end{equation}
Fix some point $q_1\in \mathcal F_s(q_0)$ and denote by 
${\mathcal G}_s(q_1)$ the integral manifold of the 
distribution ${\rm span} \{D_l: 1\leq l\leq N, l\neq s\}$, 
passing through $q_1$. Fix some vector $v\in D_s(q_1)$ such 
that $g_s(v,v)=1$. By construction for any $q\in{\mathcal 
G}_s(q_1)$ there exist a unique vector $Y_v(q)\in D_s(q)$ 
such that $d (pr_s)_ q Y_v(q)=v$. Denote by $\varepsilon_v$ 
the following function on ${\mathcal G}_s(q_1)$ 
\begin{equation}
\label{eps} 
\varepsilon_v(q)\stackrel{def}{=}{G_1}_q\bigl(Y_v(q),Y_v(q)\bigr),\quad 
q\in{\mathcal G}_s(q_1) 
\end{equation} 
It is clear that the relation (\ref{metorth}) is equivalent 
to the fact that the function $\varepsilon_v$ does not 
depend on the choice of the unit vector $v$ from 
$D_s(q_1)$. 
Then in order to obtain (\ref{metorth}) on ${\mathcal 
G}_s(q_1)$, we will put 
\begin{equation}
\label{gammamu} \gamma_s=\varepsilon_v.
\end{equation}
Let us prove that the function $\varepsilon_v$ does not 
depend on unit vector $v$ from $D_s(q_1)$. Fix some $l\neq 
s$ and some vector field $Z\in D_l$, which is unit w.r.t. 
the metric $G_1$, i.e. $G_1(Z,Z)=1$. For some $j\in I_s$, 
$i\in I_l$ take an adapted frame $(X_1,\ldots, X_n)$ such 
that 
\begin{eqnarray}
&~&\label{extension}
 X_j(q)=\varepsilon_v^{-1/2}Y_v(q),\quad  \forall q\in {\mathcal 
G}_s(q_1),\\ \label{Zeq} &~& X_i=Z.  
 \end{eqnarray} First by construction one has 
\begin{equation}
\label{interpr} 
c_{ji}^j=-\frac{1}{2}\frac{X_i(\varepsilon_v)}{\varepsilon_v} 
\end{equation}
Indeed, let $(x_1,\ldots, x_n)$ be coordinates of Corollary 
\ref{coord} and suppose that $v=\sum_{k\in I_s} 
v_k\frac{\partial}{\partial x_k}$. 
 Then by 
(\ref{extension}) on ${\mathcal G}_s(q_1)$ the fields $X_j$ 
with $j\in I_s$ have the form 
\begin{equation}
\label{Xjcoord} X_j=\varepsilon_v^{-1/2}\sum_{k\in 
I_s}v_k\frac{\partial}{\partial x_k}, 
\end{equation}
while by construction $X_i\in {\rm span} 
\Bigl(\frac{\partial}{\partial x_{\bar l}}\Bigr)_ {\bar 
l\not\in I_s}$, which together with (\ref{Xjcoord}) implies 
(\ref{interpr}). On the other hand, by (\ref{cor1}) we have 
\begin{equation}
\label{indep1} 
c_{ji}^j=\frac{1}{2}X_i\left(\frac{\lambda_s}{\lambda_l}\right) 
\left(1-\frac{\lambda_s}{\lambda_l}\right)^{-1}, 
\end{equation}
where as before $\lambda_s(q)$, $\lambda_l(q)$ are the 
eigenvalues of the transition operator $S_q$, corresponding 
to the eigenspaces $D_s(q)$ and $D_l(q)$. 
So, from (\ref{interpr}), (\ref{indep1}), (\ref{Zeq}), and 
definition of $\epsilon_v$ it follows that 
\begin{eqnarray}
&~&\label{epscor1} \frac{Z(\varepsilon_v)}{\varepsilon_v}=- 
Z\left(\frac{\lambda_s}{\lambda_l}\right) 
\left(1-\frac{\lambda_s}{\lambda_l}\right)^{-1},\\ 
&~&\label{triv}\varepsilon_v(q_1)=1 
\end{eqnarray}
The right-hand side of (\ref{epscor1}) does not depend on 
the choice of the vector $v$. 
 Hence from (\ref{epscor1})-(\ref{triv}) it 
follows that on the curve $e^{tZ}q_1$ the function 
$\epsilon_v$ does not depend on the choice of the vector 
$v$. 
%
%
Note that any point of ${\mathcal G}_s(q_1)$ can be 
connected with $q_1$ by some finite concatenation of the 
integral curves of the fields $\pm Z$, where $Z\in D_l$, 
$l\neq s$. Therefore by induction on the number of 
"switches", one gets from (\ref{epscor1}) that on the 
manifold ${\mathcal G}_s(q_1)$ the function $\varepsilon_v$ 
does not depend on the choice of the vector $v$. 
Defining $\gamma_s$, as in (\ref{gammamu}), we obtain 
(\ref{metorth}) on ${\mathcal G}_s(q_1)$ and hence on $U$, 
which completes the proof of Lemma \ref{lemma1}.$\Box$ 
\medskip 

\begin{lemma}
\label{lemma2}
There exist functions $\beta_s$ on ${\mathcal F}_s(q_0)$ such that (\ref{lambdas}) 
holds. 
\end{lemma}

{\bf Proof.} Let, as above, $(x_1,\ldots,x_n)$ be some 
coordinates from Corollary \ref{coord}. Denote by $\chi_s$ 
the following $|I_s|$-tuple: 
\begin{equation}
\label{chi} \chi_s=\{x_i\}_{i\in I_s}. 
\end{equation}
Since by construction
\begin{equation}
\label{Ispart} {\rm span}\,\{X_i\}_{\in I_s}={\rm span}\,
\{\frac{\partial}{\partial x_i}\}_{i\in I_s}=D_s, 
\end{equation} relations (\ref{cor2}) and (\ref{cor3}) are 
equivalent to the following relations respectively 
\begin{eqnarray}
&~&\label{cor2l} \forall 1\leq s\neq l\leq N,\,\, i\in I_s: 
\quad \frac{\partial}{\partial 
x_i}\left(\frac{\lambda_l^2}{\lambda_s}\right)=0,\\ 
&~&\label{cor3l}
\forall 1\leq s,l,r\leq N, \,
l\neq s,\,r\neq s,\, i\in I_s:\quad 
\frac{\partial}{\partial 
x_i}\left(\frac{\lambda_l}{\lambda_r}\right)=0. 
\end{eqnarray} 
First suppose that $N=2$. 
Then from (\ref{cor2l}) there exist functions $\bar\beta_s(\chi_s)$, $s=1,3$ such that
\begin{equation}
\label{l12}
\frac{\lambda_2^2}{\lambda_1}=\bar\beta_2(\chi_2),\quad \frac{\lambda_1^2}{\lambda_2}=\bar\beta_1(\chi_1),
\end{equation}
which easily implies (\ref{lambdas}), if we take $\beta_1=\bar\beta_1^{1/3}$, $\beta_2=\bar\beta_2^{1/3}$. 
For $N>2$ a standard analysis of conditions (\ref{cor3l}) implies
that there exist functions $\beta_s(\chi_s)$ such
that
\begin{equation}
\label{var1}
\frac{\lambda_s(q)}{\lambda_l(q)}=\frac{\beta_s(\chi_s)}{\beta_l(\chi_l)}
\end{equation}
Substituting the last relation in (\ref{cor2}) one can
obtain easily that \begin{equation} \label{var2} 
\frac{\partial}{\partial 
x_j}\left(\frac{\lambda_s(q)}{\beta_l(\chi_l)}\right)=0, 
\quad j\in I_l,\, l\neq s 
\end{equation} Using standard arguments of
"separation of variables" for the last equations, one can
easily conclude that there exist a function 
$\sigma(\chi_s)$ such that 
\begin{equation} \label{var3}
\lambda_s=\sigma(\chi_s)\prod_{l\neq s}\beta_l(\chi_l).
\end{equation}
Substituting the last equation to (\ref{var1}) we obtain
that
$$\frac{\sigma_s(\chi_s)}{\sigma_l(\chi_l)}=\frac{\beta_s^2(\chi_s)}
{\beta_l^2(\chi_l)},$$
which in turn implies that $\sigma_i=C\beta_i^2$ for some
constant $C>0$. Replacing functions $\beta_i$ by $k\beta_i$
for some constant $k>0$ one can make $C=1$. So,
\begin{equation}
\label{alpha}
\lambda_s=\beta_s(\chi_s)\prod_{l=1}^N\beta_l(\chi_l),
\end{equation}
which is equivalent to (\ref{lambdas}). $\Box$
\medskip

\begin{lemma}
If $\dim {\mathcal F}_s>1$, then $\lambda_s$ is constant on each leaf of the foliation ${\mathcal F}_s$
\end{lemma}

{\bf Proof.}
Taking the coefficients of $u_i^2$, $i\neq j$ , from (\ref{Rjexp})
and using (\ref{Rj0}), we obtain the following relation
\begin{equation}
\label{cor1ad} X_j(\alpha_i^2)=2 c_{ji}^i
(\alpha_j^2-\alpha_i^2)\quad i\neq j.
\end{equation}
Note that identity 
(\ref{cor1ad})
is stronger than identity
(\ref{cjii}):
in the first
identity we assume that the corresponding indices are
different, while in the second one
we assume that the corresponding eigenvalues are different.
 Take any pair of indices $i,j\in I_s$ such
that $i\neq j$ (by assumption $|I_s|>1$ it is
possible). Applying (\ref{cor1ad}) and using the fact that
$\alpha_i=\alpha_j=\lambda_s^{1/2}$,
 we get
 $X_j(\lambda_s)=0$ for any $j\in I_s$, which implies the statement of the lemma. 
$\Box$

\begin{remark}
\label{XYrem} The functions $\beta_s$ from relation (\ref{lambdas})
have the intrinsic meaning, because they can be expressed by 
the eigenvalues of the transition operator $S_q$
in the following way
\begin{equation}
\label{XYeq}
\beta_s\circ pr_s=\lambda_s^{^\frac{N-1}{N+1}}\Bigl(\prod_{l\neq
s}\lambda_l\Bigr)^{-\frac{2}{N+1}} \end{equation}
\end{remark}
\medskip

From the previous lemma and (\ref{lambdas}) it follows immediately the following
\begin{cor}
If $\dim {\mathcal F}_s>1$, then the function $\beta_s$ is constant.
\end{cor}

To complete the "only if" part it remains to prove relation 
(\ref{gammas}). For this, combining (\ref{gammamu}), 
(\ref{interpr}), and (\ref{indep1}), then taking into 
account (\ref{Ispart}) and (\ref{var1}), one obtains 
without difficulties 
\begin{equation}
\label{indep2}
 \forall 1\leq s\neq l\leq N,\,\,i\in I_l:\quad 
 \frac{\partial}{\partial x_i}\ln \gamma_s= 
 \frac{\partial}{\partial x_i}\ln\left 
|\frac{\beta_s(\chi_s)}{\beta_l(\chi_l)}-1\right|.
\end{equation}
 Again using standard 
"separation of variables" arguments we get from the last 
relations that there exist one-valuable functions 
$\omega_s(\chi_s)$ such that 
\begin{equation}
\label{aj} \gamma_s=\omega_s(\chi_s) \prod_{l\neq
s}\left|\frac{1}{\beta_l(\chi_l)}-\frac{1}{\beta_s(\chi_s)}
\right|. \end{equation}
 Finally note that by a change of coordinates of the type
 $\chi_s\mapsto F_s(\chi_s)$ we can make $\omega_s\equiv 1$
 for any $1\leq s\leq N$, which together with (\ref{aj}) implies 
 (\ref{gammas}). This completes the proof of the "only if" part.

Note that in the proof of the "only if" part we actually 
have used all information, which can be obtained from 
relations (\ref{Rj0}) (the only group of coefficients in 
(\ref{Rjexp}) that we did not exploit are coefficients of 
$u_iu_j$ with $i\neq j$, but the identities that they 
produce from (\ref{Rj0}) are equivalent to identities 
(\ref{cor1ad}), which was obtained by exploiting another 
group of coefficients). Therefore by Lemma \ref{prepLC} the 
conditions of the theorem are not only necessary, but also 
sufficient. 
The proof of the theorem is completed. $\Box$ \medskip

For metrics on surfaces Levi-Civita's theorem is called 
also Dini's Theorem, because Dini obtained it first in 
\cite{dini}. 


\section{The case of corank one distributions}
\setcounter{equation}{0} \indent

In the present section we investigate the problem of 
geodesic equivalence of sub-Riemannian metrics on a 
distribution $D$ of corank 1, especially, if $D$ is contact 
or quasi-contact. From the beginning we work in the 
neighborhood of regular point $q_0$, extending then the 
results to the non-regular points by the limiting process, 
when it is possible. 

Let the functions $R_j$ and $Q_{jk}$ be as in (\ref{R}) and 
(\ref{Q}) respectively. All these functions are polynomials 
on the fibers. In general, these functions depend on the 
choice of the adapted frame to the pair of the metrics 
$(G_1,G_2)$.  
\begin{defin}
\label{2divcond} We will say that the ordered pair 
$(G_1,G_2)$ of sub-Riemannian metrics on the distribution 
$D$ satisfies the second divisibility condition on an open 
set $U$, if there exist an adapted frame to the pair 
$(G_1,G_2)$ in $U$ such that for any $q\in U$ on the fiber 
$T_q^*U$ the polynomial $R_j$ is divided by the polynomial 
$Q_{j m+1}$ for any index $j$ such that $Q_{j 
m+1}\not\equiv 0$ on $T_q^*U$, $1\leq j\leq m$. 
\end{defin} 

Note that $\bar 
c_{ji}^{m+1}=\frac{1}{\alpha_i\alpha_j}c_{ji}^{m+1}$ for 
any $i,j$ such that $1\leq i\leq j$. Therefore 
\begin{equation}
\label{Qm+1} 
Q_{jm+1}=\frac{1}{\alpha_j}\sum_{i=1}^mc_{ji}^{m+1}u_i. 
\end{equation} 
\begin{prop}
\label{2divprop} Suppose that for given two sub-Riemannian 
metrics $G_1$ and $G_2$ on corank 1 distribution $D$ and 
for some open set $U$ of regular point $q_0$ there exists 
an orbital diffeomorphism of the extremal flows of these 
metrics in some open set $\mathcal B$ in $H_1\cap T^*U$, 
 $\pi(\mathcal 
B)=U$. Then the pair $(G_1, G_2)$ satisfies the second 
divisibility condition on $U$. 
\end{prop}

{\bf Proof.} Fix some index $j$, $1\leq j\leq m$, such that 
\begin{equation}
\label{QJ0}
 Q_{jm+1}\not\equiv 0.
\end{equation} 
Substituting (\ref{Philast}) into (\ref{II}) we obtain 
$$-\cfrac{\vec h_1(Q_{jm+1})R_j}{\alpha_j 
Q_{jm+1}^2{\mathcal P}^{1/2} }=\cfrac{{\rm 
polynomial}}{Q_{jm+1}{\mathcal P}^{3/2}}$$ or, 
equivalently, \begin{equation} \label{2proofeq} 
 \cfrac{{\mathcal P} \vec h_1(Q_{jm+1})
R_j}{Q_{jm+1}}={\rm polynomial}. 
\end{equation} 

Positive definite quadratic form ${\mathcal P}$ cannot be 
divided by $Q_{jm+1}$, which is linear function with real 
coefficients. Let us prove that $Q_{jm+1}$ does not divide 
$\vec h_1(Q_{jm+1})$. Assuming the converse, one can 
conclude that the coefficients of $u_j u_{m+1}$ in the 
quadratic polynomial $\vec h_1(Q_{jm+1})$ has to be equal 
to zero (because $Q_{jm+1}$ does not depend both on $u_j$ 
and on $u_{m+1}$). On the other hand, from (\ref{h1str}) 
and (\ref{Qm+1}) it is not hard to get that this 
coefficient is equal to 
$$-\frac{1}{\alpha_j}\sum_{i=1}^m\bigl(c_{ji}^{m+1}\bigr)^2.$$
Hence $c_{ji}^{m+1}=0$ for all $1\leq i\leq m$, which 
contradicts the assumption (\ref{QJ0}). So, relation 
(\ref{2proofeq}) yields that $R_j$ has to be divided by 
$Q_{jm+1}$, i.e., the second divisibility condition 
holds.$\Box$ \medskip

\begin{prop}
\label{bigcomp} Suppose that for given two sub-Riemannian 
metrics $G_1$ and $G_2$ on some $(m,m+1)$-distribution $D$ 
and for some open set $U$ there exists an orbital 
diffeomorphism of the extremal flows of these metrics in 
some open set $\mathcal B$ in $H_1\cap T^*U$, 
 $\pi(\mathcal 
B)=U$. 
Suppose also that there exists the basis $(X_1,\ldots, 
X_m)$ of $D$ adapted to the ordered pair $(G_1,G_2)$, and 
the transition operator $S_q$ has the form $S_q={\rm 
diag}\,(\alpha_1^2(q),\ldots, \alpha_m^2(q))$ in this basis 
($\alpha_i>0$). Then the following two statements hold

\begin{enumerate}
  \item If
  \begin{equation}
  \label{goodind}
  I\stackrel{def}{=}\Bigl\{j\in \{1,\ldots, m\}:  
  [X_j, D](q)\not\subset D(q) 
  \,\,\forall q\in U\Bigr\},
  \end{equation}
  then $\alpha_i=\alpha_j$ in $U$ for all $i,j\in I$;
  \item If $\alpha\stackrel{def}{=}\alpha_j$, $j\in I$, and 
  $\bar I= \Bigl\{j\in \{1,\ldots, m\}:\alpha_j=\alpha\Bigr\}$, then 
  \begin{equation}
  \label{acon1}
  \forall j\in \bar I:\quad \quad X_j(\alpha)=0
  \end{equation}
\end{enumerate}
\end{prop}

{\bf Proof.} By Proposition \ref{2divprop} for any $j\in I$ 
the polynomial $R_j$ is divided by $\alpha_j Q_{jm+1}$. But 
by (\ref{Philast}) the polynomial $\frac{R_j}{\alpha_j 
Q_{jm+1}}$ does not depend on $j\in I$ (because it is equal 
to $\sqrt{\mathcal P}\Phi_{m+1}$). In other word,
\begin{equation}
\label{rk} R_j= \Bigl(\sum_{i=1}^{m+1}r_i u_i\Bigr) 
\alpha_j Q_{jm+1}, 
\end{equation}
where coefficients $r_i$ do not depend on $j\in I$. As a 
consequence of the last identity and (\ref{Philast}) one 
has  
\begin{equation}
\label{Phir} \Phi_{m+1}=\frac{\sum_{i=1}^{m+1}r_i 
u_i}{\sqrt{\mathcal P}}. \end{equation} Using (\ref{Rjexp}) 
and (\ref{Qm+1}), one can compare the coefficients of $u_i 
u_{m+1}$, $1\leq i \leq m$ in both sides of (\ref{rk}) to 
get $$\alpha_j^2 c_{ji}^{m+1}=r_{m+1} c_{ji}^{m+1}.$$ Since 
by definition for any $j\in I$ there exist $1\leq i\leq m$ 
such that $c_{ji}^{m+1}\neq 0$, then 
\begin{equation}
\label{arm+1} \forall j\in I:\quad \alpha_j^2=r_{m+1}
\end{equation}
 In 
other words, $\alpha_j$ does not depend on $j\in I$, which 
concludes the proof of the first statement of the 
proposition. 
 
Let us prove the second statement. From (\ref{Rjexp}) and 
the fact that $\alpha_j=\alpha_i=\alpha$ for all $i\in I$ 
it follows that \begin{equation} \label{coefui} \forall 
i\in I: \Bigl({\rm the}\,\,{\rm coefficient}\,\, {\rm 
of}\,\, u_i^2\,\, {\rm in}\,\, R_j\Bigr) =-\frac{1}{2} 
X_j(\alpha^2). \end{equation} 

If $j\in \bar I\backslash I$, then $Q_{j m+1}=0$ and by 
identity (\ref{I}) we have $R_j=0$, which together with 
(\ref{coefui}) implies that $X_j(\alpha^2)=0$.

If $j\in I$, then comparing the coefficients of $u_i^2 $, 
$i\in I$, $i\neq j$ in both sides of (\ref{rk}) and using 
relations (\ref{coefui}), (\ref{Qm+1}),
we obtain \begin{equation} \label{compsq} \frac{1}{2} 
X_j(\alpha^2)=r_ic_{ij}^{m+1}. \end{equation} Substituting 
identity (\ref{Phir}) into identity (\ref{II}) with 
$s=m+1$, then using (\ref{pi}), and finally multiplying 
both sides on $\sqrt{\mathcal P}$, we get 
\begin{equation}
\label{IIr} \vec h_1\bigl(\sum_{i=1}^{m+1}r_i 
u_i\bigr)-\frac{1}{2}\Bigl(\sum_{j=1}^m\frac{X_j(\alpha_j^2)}{\alpha_j^2}
u_j\Bigr) 
\sum_{i=1}^{m+1}r_i u_i-Q_{m+1\,m+1}\sum_{i=1}^{m+1}r_i 
u_i=\sum_{k=1}^m Q_{m+1\,k}\alpha_k u_k 
\end{equation} 
Comparing the coefficients of $u_j u_{m+1}$, $j\in I$ in 
both sides of (\ref{IIr}) one can obtain without 
difficulties that $$\sum_{i=1}^m r_i 
c_{ij}^{m+1}+\frac{1}{2}X_j(\alpha^2)=0,$$ which together 
with (\ref{compsq}) implies that $\frac{n_j+1}{2} 
X_j(\alpha^2)=0$, where $n_j$ is the number of indices $i$, 
$1\leq i\leq m$ such that $c_{ij}^{m+1}\neq 0$. Therefore 
$X_j(\alpha^2)=0$ for all $j\in I$. The proof of the second 
statement is also completed. $\Box$ \medskip

As a direct consequence of Proposition \ref{subrorb} and 
the previous proposition we obtain the following 
\begin{theor}
\label{contprop}
 If two sub-Riemannian metrics $G_1$ and $G_2$,
defined on a contact distribution D, are geodesically
equivalent at some point $q_0$, then they are constantly 
proportional in some neighborhood of $q_0$. 
\end{theor}

{\bf Proof.}
First note that it is sufficient to prove this theorem for 
regular $q_0$: using the density of the set of regular 
points (Proposition \ref{denseprop}), one can extend the 
theorem to the non-regular points by passing to the limit. 
If $q_0$ is regular, then 
 there exists the basis $(X_1,\ldots, 
X_m)$ of $D$ adapted to the ordered pair $(G_1,G_2)$. Let, 
as before, the transition operator $S_q$ has the form 
$S_q={\rm diag}\,(\alpha_1^2(q),\ldots, \alpha_m^2(q))$ in 
this basis ($\alpha_i>0$). In the case of the contact 
distribution the set $I$, defined by (\ref{goodind}), 
coincides with $\{1,\ldots, m\}$. 
 Therefore, by consecutive  use of Propositions \ref{subrorb} 
  and \ref{bigcomp} we obtain that there exists the function 
  $\alpha$ such that 
  $\alpha_i=\alpha$ and $X_i(\alpha)=0$ for any $i$, $1\leq i\leq m$.
  This together with the fact that contact distribution is bracket generating implies that 
  $\alpha_i=\alpha= const$  for any $i$, $1\leq i\leq m$, which 
  concludes the proof of the theorem. $\Box$ \medskip 
%

For $(2,3)$-distributions we can extend the last result 
from contact to all bracket-generating distributions, 
because the set of points, where bracket-generating 
$(2,3)$- distributions are contact, is open and dense. 
Namely, we have the following 

\begin{cor}
\label{23cor} If two sub-Riemannian metrics $G_1$ and 
$G_2$, defined on a bracket-generating $(2,3)$-distribution 
D, are geodesically equivalent at some point $q_0$, then 
they are constantly proportional in some neighborhood of 
$q_0$. 
\end{cor}

Now consider the case of the quasi-contact distribution 
$D$. The following theorem gives the classification of all 
geodesically equivalent sub-Riemannian metrics, defined on 
such distribution: 

\begin{theor}
Suppose that $G_1$ and $G_2$ are two sub-Riemannian metrics
on the quasi-contact distribution $D$ such that
$G_2\not\equiv \mathrm{const}\, G_1$. Assume also that the
vector field $X$ is tangent to the abnormal line
distribution of $D$ and unit w.r.t. the metric $G_1$ (i.e.,
$G_{1_q}(X,X)=1$). Then the metrics $G_1$ and $G_2$ are
geodesically equivalent at the point $q_0$ if and only if
in some neighborhood $U$ of $q_0$ 
the following four conditions hold simultaneously:
\begin{enumerate}
\item If 
\begin{equation}
\label{orthD}
 D_i(q)=\{v\in D(q):G_{i_q}(v,X)=0\},\quad  i=1,2,
\end{equation}
 then 
$D_1(q)=D_2(q)$ and the distribution $D_1^2$ is codimension 
$1$ integrable distribution (here $D_1^2=D_1+[D_1,D_1]$); 
\item If $\mathcal F$ is the foliation of the integral
hypersurfaces of the distribution $D_1^2$, then the flow
$e^{tX}$ generated by the vector field $X$ preserves the
foliation ${\mathcal F}$, i.e., it maps any leaf of
$\mathcal F$ to a leaf of $\mathcal F$;
\item There exists the
one-variable function $\beta(t)$, $\beta(0)=1$, such that
 if ${\mathcal F}_0$ is the leaf of the foliation
$\mathcal F$ passing through $q_0$ and
$G_1\Bigl|_{e^{tx}{\mathcal F}_0}\Bigr.$ is the restriction 
of the metric $G_1$ to the leaf $e^{tX} {\mathcal F}_0$, 
then 
\begin{equation}
\label{lG1} G_1\Bigl|_{e^{tX} {\mathcal 
F}_0}\Bigr.=\beta(t)\Bigl((e^{-tX})^*G_1\Bigr)\Bigl|_{e^{tX} 
{\mathcal F}_0}\Bigr.; 
\end{equation}
 \item There exist two
constants $C_1>0$ and $C_2>-1$, $C_2\neq 0$, such that if
${\mathcal F}_0$ is as before and $G_2\Bigl|_{e^{tx} 
{\mathcal F}_0}\Bigr.$ is the restriction of the metric 
$G_2$ to the leaf $e^{tX} {\mathcal F}_0$, then 
\begin{equation}
\label{lG2} G_2\Bigl|_{e^{tX} {\mathcal 
F}_0}\Bigr.=\frac{C_1}{1+C_2\beta(t)}G_1\Bigl|_{e^{tX} 
{\mathcal F}_0}\Bigr., 
\end{equation}
\begin{equation}
\label{length2}\forall q\in e^{tX} {\mathcal F}_0:\quad
G_{2_q}\bigl(X(q),
X(q)\bigr)=\frac{C_1}{\bigl(1+C_2\beta(t)\bigr)^2}.
\end{equation}
\end{enumerate}
\label{quasit}
\end{theor}

Before proving Theorem \ref{quasit}, let us make some
remarks. According to this theorem for the quasi-contact
distribution $D$ the pair $(G_1,G_2)$ of constantly
non-proportional geodesically equivalent metrics at the
point $q_0$ is uniquely determined by fixing
\begin{description}
\item[~~~a)] a vector field $X$ tangent to
the abnormal line distribution of $D$;
\item[~~~b)] a hypersurface
${\mathcal F}_0$, passing through $q_0$ and transversal to 
the abnormal line distribution of $D$; 
\item[~~~c)]a sub-Riemannian
metric $\bar G$ on the contact distribution $\bar D$,
defined on the hypersurface ${\mathcal F}_0$ as follows: 
$\bar D(q)=D(q)\cap T_q {\mathcal F}_0$, $q\in {\mathcal 
F}_0$; 
\item[~~~d)] a one-variable
function $\beta(t)$ with $\beta(0)=1$;
\item[~~~e)]
two constants $C_1$, $C_2$, where $C_1>0$ $C_2>-1$, and
$C_2\neq0$.
\end{description}

The metrics $G_1$ can be uniquely recovered from the data
of a)-d). For this we extend the distribution $\bar D$ and
the metric $\bar G$ on $\bar D$
 from ${\mathcal F}_0$ to $M$ by the flow $e^{tX}$. Namely, we set
$$\forall q\in {\mathcal F}_0:\quad \bar D(e^{tX} 
q)=(e^{tX})_*\bar D(q),\quad \bar G_{e^tX q}(v,w)= \bar 
G_q\bigl((e^{-tX})_*v, (e^{-tX})_*w\bigr), \,\,v, w\in \bar 
D(e^{tX}q).$$ Then the metric $G_1$ is uniquely defined by 
the following two conditions: 
\begin{itemize}
  \item for any $q\in e^{tX} {\mathcal F}_0$ on the subspace $\bar D(q)$ the
metric $G_1$ coincides with $\bar G$ multiplied by the
factor $\beta(t)$
  \item for any $q$ the vector $X(q)$ is unit and orthogonal to
$\bar D(q)$ w.r.t. $G_1$.
\end{itemize}
In particular, it shows that the metrics on quasi-contact
distributions, admitting constantly non-proportional
geodesically equivalent metrics, are very special. The
metric $G_2$ is uniquely determined by $G_1$ and two
constants $C_1$ and $C_2$ with the properties prescribed in 
e). In other words, if the metric $G_1$ admits constantly 
non-proportional geodesically equivalent metrics, then the 
set of such metrics is two-parametric. Note also that if 
one takes $C_2=0$ in statement 4 of Theorem \ref{quasit}, 
then the metrics are constantly proportional. 

{\bf Proof of Theorem \ref{quasit}.} Let us prove the "only
if" part. Let the metrics $G_1$ and $G_2$ be geodesically
equivalent. First suppose that $q_0$ is regular point
w.r.t. the pair $(G_1,G_2)$. As before let $(X_1,\ldots,
X_m)$ be a basis of $D$ adapted to the ordered pair
$(G_1,G_2)$ and suppose that the transition operator
$S_q={\rm 
diag}\bigl(\alpha_1^2(q),\ldots,\alpha_m^2(q)\bigr)$ (where 
$\alpha_i>0$) w.r.t. the basis ($X_1,\ldots, X_m)$. 

First note that the field $X$ has to coincide with one of 
the fields $X_i$, $1\leq i\leq m$. Otherwise, the set $I$, 
defined by (\ref{goodind}), coincides with 
$\{1,\ldots,m\}$. Then by the same arguments, as in the 
proof of Theorem \ref{contprop}, we obtain that the metrics 
$G_1$ and $G_2$ are constantly proportional, which 
contradicts our assumptions. Without loss of generality, it 
can be assumed that $X=X_m$. Secondly by Proposition 
\ref{bigcomp} for any $1\leq i,j\leq m-1$ we have 
$\alpha_i=\alpha_j$. In the sequel we set $\alpha_i=\alpha$ 
for $1\leq i\leq m-1$. 

Since the field $X_m=X$ has no singularities, 
by passing to the limit one obtains that the adapted basis 
with the same properties exists also in a neighborhood of 
non-regular points w.r.t. to the pair $(G_1,G_2)$. 
Moreover, 
 $\alpha_m\neq \alpha$.
Indeed, assuming the converse we obtain from the statement 
2 of Proposition \ref{bigcomp} that the set $\bar 
I=\bigl\{j\in \{1,\ldots, m\}:\alpha_j=\alpha\bigr\}$ 
coincides with $\{1,\ldots,m\}$. But from this again by the 
same arguments, as in the proof of Theorem \ref{contprop}, 
we obtain that the metrics $G_1$ and $G_2$ are constantly 
proportional, which contradicts our assumptions. Actually, 
we have shown that for geodesically equivalent metrics 
$q_0$ is always regular: in some neighborhood of $q_0$ the 
number of distinct eigenvalues of the transition operator 
is constant and equal either to 1 (in this case the metrics 
are constantly proportional) or to 2. Besides, if $D_1$ and 
$D_2$ are as in (\ref{orthD}), then $$D_1=D_2={\rm 
span}(X_1,\ldots, X_{m-1}).$$

From (\ref{cor2}) it follows that 
$X_i\Bigl(\frac{\alpha_m^2}{\alpha}\Bigl)=0$ for all $1\leq 
i\leq m-1$, which together with (\ref{acon1}) implies 
\begin{equation}
\label{Xim} \forall 1\leq i\leq m-1:\quad X_i(\alpha_m)=0.
\end{equation}
Replacing the (\ref{acon1}) and (\ref{Xim}) in 
(\ref{cor1}), we obtain also that 
\begin{equation}
\forall 1\leq i\leq m-1:\quad \label{cmim} c_{mi}^m=0.
\end{equation} 

Let us complete the adapted basis $(X_1,\ldots, X_m)$ 
somehow to the frame $(X_1,\ldots, X_{m+1})$. 
\begin{lemma}
\label{m-1hol} 
The distribution $D_1^2=D_1+[D_1,D_1]$ is integrable.
\end{lemma}
{\bf Proof.} Using (\ref{Rjexp}) and (\ref{Qm+1}), let us 
compare the coefficients of $u_i u_m$, $1\leq i \leq m-1$ 
in both sides of (\ref{rk}), where $1\leq j\leq m-1$. As a 
result, we get easily 
$$\forall 1\leq i\neq j\leq m-1:\quad 
(\alpha^2-\alpha_m^2)c_{ji}^m=r_mc_{ji}^{m+1}+r_ic_{jm}^{m+1}.$$ 
But by construction $m\not\in I$, i.e., $c_{jm}^{m+1}=0$ 
for all $1\leq j\leq m-1$. Therefore the last relation is 
equivalent to the following identity: 
\begin{equation}
\label{D2prep} \forall 1\leq i\neq j\leq m-1:\quad 
c_{ji}^m=\frac{r_m}{\alpha^2-\alpha_m^2}c_{ji}^{m+1}. 
\end{equation} 
Hence $[X_i, X_j]\in {\rm span}\Bigl(X_1,\ldots, X_{m-1}, 
\frac{r_m}{\alpha^2-\alpha_m^2}X_m+X_{m+1}\Bigr)$ for all 
$1\leq i,j\leq m-1$ or , equivalently, \begin{equation} 
\label{D1sq} D_1^2={\rm span}\Bigl( 
D_1,\frac{r_m}{\alpha^2-\alpha_m^2}X_m+X_{m+1}\Bigr). 
\end{equation}
To prove the lemma it is sufficient to prove that 
\begin{equation}
\label{XiXm} \forall 1\leq i\leq m-1:\quad 
\Bigl[X_i,\frac{r_m}{\alpha^2-\alpha_m^2}X_m+X_{m+1}\Bigr]\in 
{\rm span}\Bigl( 
D_1,\frac{r_m}{\alpha^2-\alpha_m^2}X_m+X_{m+1}\Bigr). 
\end{equation}
Using (\ref{acon1}) and (\ref{Xim}), it is easy to show 
that (\ref{XiXm}) is equivalent to the following 
identity
\begin{equation}
\label{idint}
 X_i(r_m)+r_m c_{mi}^m+c_{m+1 
i}^m(\alpha^2-\alpha_m^2)-r_m c_{m+1\, i}^{m+1}=0
\end{equation}

Let us prove identity (\ref{idint}). First note that from 
(\ref{acon1}) and (\ref{compsq}) it follows easily that 
$r_i=0$ for $1\leq i\leq m-1$ (here we use also the fact 
that for given $i$, $1\leq i\leq m-1$, there exist $j$, 
$1\leq j\leq m-1$, such that $c_{ij}^{m+1}\neq 0$). From 
this and (\ref{acon1}) it follows that the identity 
(\ref{IIr}) can be rewritten in the following form: 
\begin{equation}
\label{IIrm} \vec h_1\bigl(\sum_{i=m}^{m+1}r_i 
u_i\bigr)-\frac{1}{2}\frac{X_m(\alpha_m^2)}{\alpha_m^2} u_m 
\sum_{i=m}^{m+1}r_i u_i-Q_{m+1\,m+1}\sum_{i=m}^{m+1}r_i 
u_i=\sum_{k=1}^m Q_{m+1\,k}\alpha_k u_k 
\end{equation}  
Comparing the coefficients of $u_i u_m$, $1\leq i\leq m-1$ 
in both sides of (\ref{IIrm}) by use of (\ref{h1str}) and 
(\ref{Q}) it is not difficult to obtain \begin{equation} 
\label{idint1} 
 X_i(r_m)+r_m c_{mi}^m+r_{m+1}(c_{m+1 
i}^m+c_{m+1 m}^i)-r_m \bar c_{m+1\, i}^{m+1}\alpha =(\bar 
c_{m+1 m}^i+\bar c_{m+1 i}^m)\alpha \alpha_m 
\end{equation}
From (\ref{relc2}) and (\ref{alphag}) it follows that $\bar 
c_{m+1\, i}^{m+1}=\frac{1}{\alpha}c_{m+1\, i}^{m+1}$, $\bar 
c_{m+1 m}^i =\frac{\alpha}{\alpha_m}c_{m+1 m}^i$, and $\bar 
c_{m+1 i}^m=\frac{\alpha_m}{\alpha}c_{m+1 i}^m$, while by 
(\ref{arm+1}) we have $r_{m+1}=\alpha^2$. Substituting all 
this to (\ref{idint1}) we get (\ref{idint}), which 
completes the proof of the lemma. $\Box$ \medskip

\begin{lemma}
\label{preflow} If $\mathcal F$ is the foliation of the  
integral hypersurfaces of the distribution $D_1^2$, then 
the flow $e^{tX}$ generated by the vector field $X$ 
preserves the foliation ${\mathcal F}$. 
\end{lemma}

{\bf Proof.} From the previous lemma it follows that in 
some neighborhood $U$ of $q_0$ there exist coordinates 
$(x_1,\ldots, x_{m+1})$ such that the leaves of $\mathcal 
F$ are $\{x_m=const\}$ and $X_m=\nu 
\frac{\partial}{\partial x_m}$ for some function $\nu$. By 
construction, all vector fields $X_i$ with $1\leq i\leq 
m-1$ belong to ${\rm span}\bigl(\frac{\partial}{\partial 
x_1},\ldots,\frac{\partial}{\partial 
x_{m-1}},\frac{\partial}{\partial x_{m+1}}\bigr)$. 
Therefore 
$c_{mi}^m=\frac{X_i(\nu)}{\nu}$ 
for all $1\leq i\leq m-1$, which together with (\ref{cmim}) 
implies that 
$ X_i(\nu)=0$ for all $1\leq i\leq m-1$.
 Then $\nu$ is constant on each leaf of ${\mathcal F}$, which 
 equivalent to the statement 
 of the lemma. $\Box$ \medskip
 
 \begin{lemma}
 \label{cond3l}
 Relation (\ref{lG1}) holds for some one-variable function 
 $\beta(t)$.
 \end{lemma}
 
 {\bf Proof.} Actually the proof of this lemma is very 
 similar to the proof of Lemma \ref{lemma1}. 
 Since the vector field $X=X_m$ satisfies $[X,D]\subset D$, the flow 
 $e^{tX}$ preserves the distribution $D$. This and the 
 previous lemma implies that  $e^{tX}$ preserves also the 
 distribution $D_1$ (note that by the previous lemma 
 $D_1(q)=D(q)\cap T_q{\mathcal F}(q)$, where ${\mathcal 
 F}(q)$ is the leaf of the foliation ${\mathcal F}$, 
 passing through the point $q$).
 
 Fix some point $q_1\in {\mathcal F}_0$. Denote by $L_{q_1}$ the abnormal 
 extremal trajectory passing through $q_1$.  
 Fix some vector $v\in D_s(q_1)$ such 
that ${G_1}_{q_1}(v,v)=1$. By construction for any point 
$q\in L_{q_1}$ such that $q=e^{tX}q_1$ there exist a unique 
vector $Y_v(q)\in D_1(q)$ such that $d (e^{-tX})_ q 
Y_v(q)=v$. Denote by $\varepsilon_v$ the following function 
on the curve $L_{q_1}$.
\begin{equation}
\label{eps1} 
\varepsilon_v(q)\stackrel{def}{=}{G_1}_q\bigl(Y_v(q),Y_v(q)\bigr),\quad 
q\in L_{q_1} 
\end{equation} 
By the same arguments as in the proof of Lemma 
\ref{lemma1}, we obtain that the function $\varepsilon_v$ 
does not depend on the choice of the unit vector $v$ from 
$D_1(q_1)$. It implies that for any $q$ in some 
neighborhood $U$ of $q_0$ (here any coordinate neighborhood 
from the proof of the previous lemma can be taken as $U$) 
there is $\beta(q)$ such that if $q=e^{tX}q_1$, where 
$q_1\in {\mathcal F}_0$, then 
\begin{equation}
\label{lG11}
 G_1\Bigl|_{e^{tX} 
{\mathcal 
F}_0}\Bigr.=\beta\Bigl((e^{-tX})^*G_1\Bigr)\Bigl|_{e^{tX} 
{\mathcal F}_0}\Bigr. 
\end{equation}
Besides, similarly to (\ref{epscor1})-(\ref{triv}), we have
\begin{eqnarray}
&~&\label{betcor1} \frac{X(\beta)}{\beta}=- 
X\left(\frac{\alpha^2}{\alpha_m^2}\right) 
\left(1-\frac{\alpha^2}{\alpha_m^2}\right)^{-1},\\ 
&~&\label{betriv}\beta\Bigl|_{{\mathcal F}_0}\Bigr.=1.
\end{eqnarray}
 
Finally by (\ref{acon1}) and (\ref{Xim}) the functions 
$\alpha$ and $\alpha_m$ are constant on each leaf of the 
foliation ${\mathcal F}$. Therefore 
(\ref{betcor1})-(\ref{betriv}) implies that the function 
$\beta$ is constant on each leaf of the foliation 
${\mathcal F}$ too. This fact together with (\ref{lG11}) 
implies (\ref{lG1}).$\Box$\medskip

In order to complete the proof of the "only if" part it 
remains to prove identities (\ref{lG2}) and 
(\ref{length2}). By (\ref{barXdef}) and statement 1 of 
Proposition \ref{bigcomp} 

\begin{equation}
\label{lG2p} \forall q\in e^{tX} {\mathcal F}_0:\quad 
{G_2}_q =\alpha^2(q){G_1}_q 
\end{equation}
\begin{equation}
\label{length2p}\forall q\in e^{tX} {\mathcal F}_0:\quad 
G_{2_q}\bigl(X(q), X(q)\bigr)=\alpha_m^2(q). 
\end{equation} 

So, it remains to find the functions $\alpha$ and 
$\alpha_m$. As was mentioned in the proof of the previous 
lemma, the functions $\alpha$ and $\alpha_m$ are constant 
on each leaf of the foliation ${\mathcal F}$. Besides, by 
(\ref{cor2}) we have $X_m(\frac{\alpha^2}{\alpha_m})=0.$ 
So, 
\begin{equation}
\label{fraconst}
 \frac{\alpha^2}{\alpha_m}\equiv C,
\end{equation}
where $C$ is constant. Then from (\ref{cor1}) it follows 
that for any $j$, $1\leq j\leq m-1$
\begin{equation}
\label{cor1qua} 
X_m\left(\frac{C}{\alpha_m}\right)=2c_{jm}^j 
\left(1-\frac{C}{\alpha_m}\right) 
\end{equation} 
By Lemmas \ref{m-1hol}-\ref{cond3l} we can choose the 
coordinates $(y_1,\ldots, y_m, t)$ in a neighborhood of 
$q_0$ such that $q_0=(0,\ldots,0)$ and 
\begin{eqnarray}
&~&\label{dt} X_m=\frac{\partial}{\partial t}; \\ &~& 
\label{Xnu} 
X_j=\beta(t)^{-1/2}\sum_{k=1}^m\nu_{jk}\frac{\partial}{\partial 
y_k},\quad 1\leq j\leq m-1. 
\end{eqnarray}
As in (\ref{interpr}) this yields that 
$$c_{jm}^j=-\frac{1}{2}\frac{d}{dt}\ln \beta(t).$$ 
Substituting the last formula in (\ref{cor1qua}), one can 
obtain without difficulties that  
\begin{equation} \label{alpham} \alpha_m =\frac{C}
{1+C_2\beta(t)}
\end{equation}
for some constant $C_2$, $C_2>-1$, $C_2\neq 0$. Then by 
(\ref{fraconst}) 
\begin{equation} \label{alphaC} \alpha^2 =C\alpha_m=\frac{C^2}
{1+C_2\beta(t)} 
\end{equation}
Setting $C_1=C^2$ and substituting (\ref{alphaC}) and 
(\ref{alpham}) into (\ref{lG2p}) and (\ref{length2p}), we 
get (\ref{lG2}) and (\ref{length2}). The proof of the "only 
if" part of the theorem is completed. 

Note that in the proof of the "only if " part we actually 
have used all information, contained in (\ref{rk}), which 
is equivalent to (\ref{I}). Also, it can be shown by direct 
check that if all conditions 1-4 of Theorem \ref{quasit} 
hold then the identity (\ref{IIrm}) holds too (but this 
identity is equivalent to (\ref{II})). From this, Lemma 
\ref{lastcomp}, and Proposition \ref{sufdif} it follows 
that conditions 1-4 of the theorem are also sufficient for 
the geodesic equivalence of our metrics at $q_0$. $\Box$

\section{
The case of Riemannian metrics on a surface near
non-regular isolated point} \setcounter{equation}{0} 
\indent 

 In the present
section for the Riemannian metrics on a surface we obtain 
the classification of geodesically equivalent pairs at 
non-regular point (the point of bifurcation of the 
eigenvalues of the transition operator). Namely, we 
consider the case when two Riemannian metrics on a surface 
are proportional in an isolated point. Since the set of all 
$2\times 2$ symmetric matrices with the equal eigenvalues 
has codimension 2 in the set of all $2\times 2$ symmetric 
matrices, we have that for generic pair of Riemannian 
metrics on a surface the set of points of their 
proportionality consists of isolated points. Therefore its 
is natural to consider the case when two Riemannian metrics 
on a surface are proportional in an isolated point. It 
turns out that 
 Dini's
Theorem (i.e., Levi-Civita's theorem in the case of a  
surface) can be naturally extended to this case. 

First let us formulate Dini's Theorem 
in the case of non-proportional metrics and analyze its 
additional features. 

\begin{theor}
(Dini's Theorem) 
\label{dinit1} Suppose that two Riemannian metrics $G_1$ 
and $G_2$ on a surface are non-proportional at some point 
$q_0$. Then they are geodesically equivalent at $q_0$ if 
and only if in some neighborhood of $q_0$, there exist 
coordinates $(x_1, x_2)$, $q_0=(x_1^0,x_2^0)$, and 
one-variable functions $\beta_1(x_1)$ and $\beta_2(x_2)$ 
($\beta_1(x_1^0)<\beta_2(x_2^0)$) such that in this 
coordinates 
\begin{equation}
\label{met1d}
||\cdot||_1^2=\left(\frac{1}{\beta_1(x_1)}-\frac{1}{\beta_2(x_2)}\right)
(dx_1^2+dx_2^2),
\end{equation}
\begin{equation}
\label{met2d}
||\cdot||_2^2=\beta_1(x_1)\beta_2(x_2)\left(\frac{1}{\beta_1(x_1)}-
\frac{1}{\beta_2(x_2)}\right) \left(\beta_1(x_1) dx_1^2+ 
\beta_2(x_2)dx_2^2 \right), 
\end{equation}
where $||v||_i^2=G_i(v,v)$, $i=1,2$. 
\end{theor}
%
%

The coordinates, appearing in Theorem \ref{dinit1}, will be 
called {\it Dini's coordinates} of the ordered pair of 
Riemannian metrics $(G_1, G_2)$. The following lemma will 
be useful in the sequel 

\begin{lemma}
\label{dinitrans} If $(x_1,x_2)$ and $(\bar x_1, \bar x_2)$ 
are two Dini's coordinates of the ordered pair of 
Riemannian metrics $(G_1,G_2)$ on the same neighborhood 
$U$, then $\bar x_i=\pm x_i+c_i$
some constants $c_i$, $i=1,2$. 
\end{lemma}

{\bf Proof.} From Corollary \ref{coord} and the fact that 
in Theorem \ref{dinit} we assume that 
$\beta_1(x_1^0)<\beta_2(x_2^0)$ it follows that the 
coordinate net of all Dini's coordinates on $U$ coincide: 
$D_1=\{dx_2=0\}=\{d\bar x_2=0\}$ is the line distribution 
of the eigenvectors of the transition operator, 
corresponding to its smallest eigenvalue, while 
$D_2=\{dx_1=0\}=\{d\bar x_1=0\}$ is the line distribution 
of the eigenvectors of the transition operator, 
corresponding to its biggest eigenvalue. Hence 
the transition function between the coordinates has a form 
$x_i=\psi_i(\bar x_i)$, $i=1,\ldots n$. Then the first 
metric is written in the coordinates $(\bar x_1,\ldots,\bar 
x_n)$ as follows: $$||\cdot||_1^2=
\left(\frac{1}{\beta_1\bigl(\psi(\bar 
x_1)\bigr)}-\frac{1}{\beta_2^2\bigl(\psi (\bar x_2)\bigr)} 
\right)\sum_{j=1}^2(\psi_i'(\bar x_j))^2 (d\,\bar x_j)^2.$$ 
By Remark \ref{XYrem} the coefficients of $d x_j^2$ in 
(\ref{met1d}) do not depend on the choice of the Dini 
coordinates. Therefore $(\psi_i'(\bar x_j))^2\equiv 1$, 
which implies the statement of the Lemma. $\Box$ 


 Recall that a Riemannian 
metric on a surface defines the canonical conformal 
structure: In a neighborhood of any point there is a 
coordinate system in which the Riemannian metric has the 
form \begin{equation} \label{isoterm} 
||\cdot||^2=\mu(x,y)(d\,x^2+d\,y^2). 
\end{equation} Such coordinates are called {\it isothermal}
(see, for example, \cite{spr} or \cite{eisen}). The
transition function from one isothermal coordinates to some 
other is conformal mapping, up to the orientation, so the 
set of all charts with isothermal coordinates defines the 
conformal structure. 
Note that by (\ref{met1d}) all Dini's coordinates are 
isothermal w.r.t. the first metric $G_1$.

Now suppose that the Riemannian metrics $G_1$ and $G_2$ are 
proportional at some isolated point $q_0$ and geodesic 
equivalent in a neighborhood of this point. Choose in a 
neighborhood $B$ of $q_0$ some isothermal coordinates 
$(x,y)$ w.r.t. the first metric $G_1$. 
Also, we can assume that the metrics are geodesic
equivalent in $B$ (otherwise we can take a smaller
neighborhood).
By above for any $q\in B$ in a neighborhood $B_q$ of $q$
there exist Dini's coordinates $(u,v)$ of the ordered pair
$(G_1,G_2)$. We also can take them such that they define 
the same orientation as $(x,y)$. The pair $(B_q, 
u(x,y)+iv(x,y)$ is a function element of an analytic 
function. Taking one of such function elements and using 
the standard procedure of the analytic continuation, we get 
the analytic function $F$ in the punctured neighborhood 
$B\backslash q_0$ such that each of its function elements 
defines the transition function from the chosen isothermal 
coordinates $(x,y)$ to  
Dini's coordinates of the ordered pair $(G_1,G_2)$ in the 
neighborhood of this function element. The function $F$ 
will be called a {\it Dini transition function} of the 
ordered pair of geodesic equivalent Riemannian metrics from 
the given isothermal coordinates $(x,y)$. The following 
theorem gives the characterization of Dini's transition 
functions at an isolated point of the proportionality of 
the metrics: 

\begin{theor}
\label{bift} If $F(z)$ is some Dini transition function of 
the ordered pair of Riemannian metrics, which are 
proportional at an isolated point $q_0$ and geodesic 
equivalent in a neighborhood of this point, then the 
function $(F')^2$ has a pole of order 1 or 2 at $q_0$. 
Besides, if $(F')^2$ has a pole of order 2 at $q_0$, then 
the principle negative coefficient in its Laurent expansion 
at $q_0$ has to be real. 
\end{theor}

{\bf Proof.} First note that the function $(F')^2$ is an
one-valued function on some punctured neighborhood
$B\backslash q_0$ of $q_0$. Indeed, by Lemma
\ref{dinitrans} the function elements $(V, F_1)$ and $(V,
F_2)$ of $F$ (with the common neighborhood of definition)
satisfy
\begin{equation}
\label{br} F_1(z)\equiv\pm F_2(z)+c,\quad z=x+iy,
\end{equation}
 where $c$
is some complex constant. This implies that $(F_1')^2\equiv
(F_2')^2$.

Now let us prove that $(F')^2$ has a pole at $q_0$.
Indeed, suppose that in the original coordinates the metric
$G_1$ satisfies (\ref{isoterm}) with some function $\mu$. 
Writing the metric $G_1$ in Dini's coordinates, we obtain 
that 
\begin{equation}
\label{lam} \mu=\Bigl(\frac{1}{\beta_1}-
 \frac{1}{\beta_2}\Bigr)|F'|^2,
 \end{equation}
where the functions $\beta_i$ are as in Theorem
\ref{dinit1}. The functions $\beta_i$ are expressed by the
eigenvalues $\lambda_j$ of the transition operator as in
(\ref{XYeq}) with $n=2$. The condition of the 
proportionality of the metrics at $q_0$ implies that 
\begin{equation}
\label{betq0}
\beta_1(q_0)=\beta_2(q_0)
\end{equation}
(because $\lambda_1(q_0)=\lambda_2(q_0)$). From this, 
(\ref{lam}) and the fact that the function $\mu$ has no 
singularity at $q_0$ it follows that 
$\lim_{z\rightarrow q_0}|F'(z)|^2=\infty$, i.e. $(F')^2$
has a pole at $q_0$.

Although the function $F$ is in general multiple-valued, by
(\ref{br}) the families of the level sets of the function
${\rm Re}\, F$ (the function ${\rm Im}\, F$) for all its
branches coincide. By construction the function $\beta_1$
is constant on the level set of ${\rm Re}\, F$, while the
function $\beta_2$ is constant on the level set of ${\rm
Im}\, F$. Using this fact it is not difficult to prove that 
the order of pole of $(F')^2$ at $q_0$ is not greater than 
$2$. Assuming the converse, we obtain that the function 
$F^2$ also has a pole at $q_0$. So, $F^2$ maps a puncture 
neighborhood of $q_0$ onto the neighborhood of infinity and 
also sends the point $q_0$ to $\infty$. But any level set 
of ${\rm Re}\, F$ is the preimage w.r.t. the mapping $F^2$ 
of some parabola of the type $u=c^2-\frac{v^2}{4c^2}$ on 
the plane $w$, where $w=F^2(z)$, $w=u+iv$. Such parabolas 
have $\infty$ as an accumulation point. Hence $q_0$ is the 
accumulation point of any level set of the function ${\rm 
Re}\, F$. This together with the fact that $\beta_1$ is 
constant on the level set of ${\rm Re}\, F$ and continuous 
at $q_0$ implies that $\beta_1$ is identically equal to 
some constant $C_1$ in a neighborhood of $q_0$. In the same 
way, $\beta_2$ is identically equal to some constant $C_2$ 
there. Moreover, by (\ref{betq0}) $C_1=C_2$. But it means 
that our metrics are proportional in the neighborhood of 
$q_0$, which contradicts our assumptions. So, the order of 
pole of $(F')^2$ at $q_0$ is not greater than $2$. 

To complete the proof of the theorem it remains to show
that if $(F')^2$ has a pole of order 2 at $q_0$, then the
principle negative coefficient in its Laurent expansion at
$q_0$ is real. Indeed, if $(F')^2$ has a pole of order 2
with the principle negative coefficient $a$ in its Laurent
expansion at $q_0$, then $F$ has the logarithmic
singularity at $q_0$ with coefficient $\sqrt {a}$ near the
logarithm. In this case after the appropriate change of
independent variable $z$ in a neighborhood of $q_0$ (i.e.,
conformal change of coordinates in a neighborhood of $q_0$)
one can get $F(z)=\sqrt{a} \log z$, $q_0=0$. But if $a$ is
not real, then $\sqrt{a}$ is neither real nor pure
imaginary. In this case all level sets of both ${\rm Re} F$
and ${\rm Im} F$ are spirals having $q_0$, as an 
accumulation point. As above, it implies that functions 
$\beta_1$ and $\beta_2$ are equal to the same constant, 
which is impossible. The proof of the theorem is completed. 
$\Box$ 
\medskip

According to the previous theorem only the following two
situation are possible at an isolated point of
proportionality of two metrics:

{\bf 1)} $(F')^2$ has a simple pole at $q_0$
$\Leftrightarrow$ $F(z)=\sqrt{G(z)}$ for some analytic
function $G(z)$, having at $q_0$ zero of order 1 (i.e., $F$
has the "square root" singularity at $q_0$). In this case
after the appropriate change of independent variable $z$ in
a neighborhood of $q_0$ one can get $F(z)=\sqrt z$,
$q_0=0$;

{\bf 2)} $F'$ has a simple pole at $q_0$ with real or pure
imaginary residue at $q_0$ $\Leftrightarrow$ $F(z)$ has the
logarithmic singularity at $q_0$ with real or pure
imaginary coefficient $b$ near logarithm. In this case
after the appropriate change of independent variable $z$ in
a neighborhood of $q_0$ one can get $F(z)=b \log z$,
$q_0=0$, where $b$ is real or pure imaginary constant. If
$b$ is real then the level sets of ${\rm Im}\, F(z)$ are
the rays, starting at $0$. Hence by the same argument as in
the proof of the previous theorem we can conclude that the
function $\beta_2$ is constant. Besides, the function
$\beta_1$ in this case depends only on $|z|$ (here
$\beta_i$ as in Theorem \ref {dinit}). In the same way, if
$b$ is pure imaginary, then $\beta_1$ is constant and
$\beta_2$ depends only on $|z|$.

Using 1), 2) and Theorem \ref{dinit}, we obtain without
difficulties the following analog of Dini's Theorem:

\begin{theor}
\label{gendini} (Generalization of Dini's Theorem to the
case of an isolated non-regular point) Two Riemannian 
metrics $G_1$ and $G_2$ on a surface $M$, which are 
proportional in an isolated point $q_0$, are geodesic 
equivalent in a neighborhood of this point if and only if 
one of the following two conditions holds: 
\begin{enumerate}
\begin{item}
In a neighborhood of $q_0$, there exist coordinates $(x,
y)$, $q_0=(0,0)$ and two one-variable smooth functions $U$
and $V$, satisfying $0<U(u)< V(0)=U(0)<V(v)$ for all
positive $u$ and $v$, $U'(0)=-V'(0)$, and $V'(0)>0$, such
that in the punctured neighborhood of $q_0$ the metrics 
$G_1$ and $G_2$ satisfy 
\begin{equation}
\label{met1ds} ||\cdot||_1^2=\left(\frac{1}{U\Bigl(r
\cos^2\frac{\theta}{2}\Bigr)}-\frac{1}{V\Bigl(r\sin^2\frac{\theta}{2}\Bigr)
} \right)\frac{1}{4 r}(dr^2+r^2 d\theta^2),
\end{equation}
\begin{equation}
\label{met2ds} ||\cdot||_2^2=
\frac{S}{8r}\bigl( (A-S \cos\theta)dr^2-2S r \sin\theta dr
d\theta+(A+S \cos\theta)r^2 d\theta^2\bigr),
\end{equation}
where $$
A=U\Bigl(r \cos^2\frac{\theta}{2}\Bigr)+V\Bigl(r
\sin^2\frac{\theta}{2}\Bigr),\quad S=V\Bigl(r
\sin^2\frac{\theta}{2}\Bigr)-U\Bigl(r
\cos^2\frac{\theta}{2}\Bigr),
$$ and $(r,\theta)$ are the corresponding polar
coordinates;
\end{item}

\begin{item}
In a neighborhood of $q_0$, there exist coordinates $(x,
y)$, $q_0=(0,0)$, positive constants $a$, $C$, and an
one-variable smooth functions $R(r)$, satisfying $R(r)\neq
R(0)$ for $r>0$, $R(0)=C$, $R'(0)=0$, and $R''(0)\neq 0$,
such that in the punctured neighborhood of $q_0$ the
metrics $G_1$ and $G_2$ satisfy
\begin{equation}
\label{met1dl} 
||\cdot||_1^2=\left|\frac{1}{C}-\frac{1}{R(r) 
}\right|\frac{a}{r^2}(dr^2+r^2 d\theta^2), 
\end{equation}
\begin{equation}
\label{met2dl} ||\cdot||_2^2=\frac{a C R(r)}{r^2}
\left|\frac{1}{C}-\frac{1}{R(r)}\right|\bigl( R( r)dr^2+C
r^2 d\theta^2\bigr),
\end{equation}
where $(r,\theta)$ are the corresponding polar coordinates.
\end{item}
\end{enumerate}

\end{theor}

\begin{remark}
The conditions on the functions $U$ and $V$ in the case 1 
and on $R$ in the case 2 follows easily from the fact that 
the metrics are positive definite and nonsingular at $q_0$. 
\end{remark}
\begin{remark}
\label{bifcurve} Using the standard arguments of Complex
Analysis, one can show that for the pair of geodesically 
equivalent metrics the set of non-regular points cannot be 
a rectifiable curve $\Gamma$: in this case one can 
construct an one-valued Dini transition conformal function 
out of $\Gamma$ which goes to infinity, when one tends to 
$\Gamma$. Then by Morera Theorem the function $1/F$ is 
analytic and equal to zero on $\Gamma$ and so everywhere, 
which is impossible. 
\end{remark}

{\it Acknowledgments.} I would like to thank professor A. 
 Agrachev for proposing me to work on this problem and for his 
 constant attention to this work.

\end{document}